\documentclass[tikz,12pt]{article}

\usepackage{epsfig}
\usepackage{latexsym}
\usepackage{caption}
\usepackage{amsfonts}
\usepackage{amssymb}
\usepackage{mathrsfs}
\usepackage{amsmath}
\usepackage{enumerate}
\usepackage{enumitem}
\usepackage{graphicx}
\usepackage{MnSymbol}
\usepackage{float}
\usepackage{pict2e}
\usepackage{graphics} 
\usepackage{pgf,tikz}

\usepackage{cite}
\usepackage{soul}                                                         

\usepackage{authblk}

\usepackage{verbatim}
\usepackage{subcaption}

\usepackage{soul,color}

\soulregister\cite7
\soulregister\ref7
\soulregister\pageref7

\definecolor{lightgreen}{RGB}{108,222,60}

\usepackage[colorinlistoftodos,prependcaption,textsize=scriptsize, textwidth=30mm,]{todonotes}

\usepackage{mathtools}

\usepackage{cite}
\usepackage{multirow} 
\usepackage{multicol} 
\usepackage{arydshln}
\usepackage{booktabs}
\usepackage{tabu}
\usepackage{ulem}
\usepackage{pifont}
\newcommand{\xmark}{\text{\ding{55}}}
\newcommand{\dmark}{\text{\ding{53}}}

\usetikzlibrary{calc}
\usetikzlibrary{arrows}

\makeatletter

\usepackage{xcolor}

\def \red {\textcolor{red} }

\usepackage[pdfauthor={derajan},pdftitle={How to do this},pdfstartview=XYZ,bookmarks=true,
colorlinks=true,linkcolor=blue,urlcolor=blue,citecolor=blue,pdftex,bookmarks=true,linktocpage=true,  hyperindex=true]{hyperref}

\makeatother

\newtheorem{theorem}{Theorem}
\newtheorem{observation}[theorem]{Observation}
\newtheorem{corollary}[theorem]{Corollary}

\newtheorem{definition}[theorem]{Definition}

\newtheorem{claim}{Claim}

\newtheorem{lemma}[theorem]{Lemma}

\makeatletter
\newcommand\figcaption{\def\@captype{figure}\caption}
\newcommand\tabcaption{\def\@captype{table}\caption}
\makeatother

\newenvironment{proof}{\noindent {\bf
		Proof.}}{\rule{3mm}{3mm}\par\medskip}

\newcommand{\fr}{\mathcal{R}}

\newcommand{\fp}{\mathcal{P}}

\newcommand{\qed}{\null\nobreak\rule{0.75em}{0.7em}}

\begin{document}

	\title{Decomposition of triangle-free planar graphs}
	
	\author[1]{Rongxing Xu \thanks{E-mail: xurongxing@ustc.edu.cn. Supported by Anhui Initiative in Quantum Information Technologies grant  AHY150200. Supported also by NSFC 11871439.}}
	\author[2]{Xuding Zhu \thanks{E-mail: xdzhu@zjnu.edu.cn. Grant Numbers: NSFC 11971438, U20A2068. ZJNSFC   LD19A010001.}}
	
	\affil[1]{\small School of Mathematical Sciences, University of Science and Technology of China, Hefei, Anhui, 230026, China}
	
	\affil[2]{\small Department of Mathematics, Zhejiang Normal University, Jinhua, Zhejiang, 321000, China.}

	\maketitle
	
	\begin{abstract}
		A decomposition  of a graph  $G$ is a family of subgraphs of $G$ whose edge sets form a partition of $E(G)$. In this paper, we prove that every triangle-free planar graph $G$ can be decomposed into a $2$-degenerate graph and a matching. Consequently, every triangle-free planar graph $G$ has  a matching $M$ such that $G-M$ is online  3-DP-colorable. This strengthens an earlier result in [R. \v{S}krekovski,  {\em A Gr\"{o}tzsch-Type Theorem for List Colourings with Impropriety One}, Combin. Prob. Comput. 8 (1999), 493-507] that every triangle-free planar graph is $1$-defective $3$-choosable.
	\end{abstract}
	
	\section{Introduction}\label{Int}
	
	A \emph{decomposition} of  a graph $G$ is a family of   subgraphs $H_1, \ldots, H_t$ such that each edge of $G$ is an edge of   exactly one of the subgraphs.  The problem of decomposing a  graph into   subgraphs with simpler structure is one of the central topics in graph theory. 
	One problem is to decompose a   graph into subgraphs of bounded degeneracy.
	The Nash-Williams Arboricity Theorem \cite{Nash-1964} (see also  \cite{Nash-1961-1} and \cite{Tutte1961}) gives a necessary and sufficient condition for a graph to be decomposable into a certain number of forests (i.e., $1$-degenerate graphs). In particular, it implies that 
	every planar graph can be decomposed into three $1$-degenerate graphs. A result by Schnyder \cite{Schnyder1990} ensures that a planar graph can be decomposed into a $1$-degenerate graph and a $2$-degenerate graph. 
	
	Another problem is to decompose a graph into a graph with bounded degeneracy (or bounded chromatic number, or bounded choice number, or bounded paint number, or bounded Alon-Tarsi number) and a graph with bounded maximum degree. This problem is related to defective coloring of graphs. A {\em $d$-defective coloring} of a graph $G$ colors its vertices in such a way that each vertex $v$ has at most $d$ neighbors that are colored the same color as $v$. Thus
	a graph $G$ is $d$-defective $k$-colorable if and only if $G$ can be decomposed into a graph of  chromatic number at most $k$ 
	and a graph of maximum degree at most $d$.
	However, $d$-defective $k$-choosable is weaker than   decomposable into a graph of  choice number at most $k$ 
	and a graph of maximum degree at most $d$.
	The latter statement asserts that $G$ has a subgraph $H$  with $\Delta(H) \le d$ and $G-E(H)$ is $L$-colorable for every $k$-list assignment $L$. The former statement asserts that
	for every $k$-list assignment $L$ of $G$, there is a subgraph $H$ of $G$ (which may depend on $L$) with $\Delta(H) \le d$ and $G-E(H)$ is $L$-colorable. 
	It was proved independently by \v{S}krekovski \cite{Riste1999-1} and Eaton and Hull \cite{ET1999} that every planar graph is $2$-defective $3$-choosable. However, it was shown in \cite{CCKPSZ2020} that there are planar graphs that cannot be decomposed into a $3$-choosable graph and a graph of maximum degree at most $3$. It was shown in \cite{CCKPSZ2020} that every planar graph can be decomposed into a $2$-degenerate (and hence $3$-choosable) graph and a graph of maximum degree at most $6$. It remains an open problem whether $6$ can be reduced to $5$ or $4$.

	In this paper, we are interested in decomposing  triangle-free planar graphs. It was proved by \v{S}krekovski \cite{Riste1999-2} that every triangle-free planar graph is $1$-defective $3$-choosable.  A natural question is whether   every triangle-free planar graph can be decomposed into a matching and a $3$-choosable graph?  This paper shows that the answer is yes. Indeed, we shall prove a stronger result.  
	By a {\em $(d,k)$-decomposition} of a graph $G$, we mean a  pair $(H_1, H_2)$, where $H_1$ is $d$-degenerate and $H_2$ has maximum degree at most $k$. We prove the following result:

	\begin{theorem}
		\label{main-thm0}
		Every triangle-free planar graph has a $(2,1)$-decomposition.
	\end{theorem}  
	
	The concepts of  {online list coloring}, DP-coloring  and online DP-coloring  of a graph are variations of list coloring. We refer readers to \cite{ DP2018,KKLZ,Schauz2009, Zhu2009}
	for the definitions and properties of these concepts. Here we just mention that online DP-$k$-colorable   implies   DP-$k$-colorable as well as  online $k$-choosable, and each of DP-$k$-colorable  and online $k$-choosable implies   $k$-choosable. On the other hand,   $k$-degenerate implies online $(k+1)$-DP-colorable.  Hence we have the following corollary, which strengthens the above mentioned result in  \cite{Riste1999-2}.
	 
	\begin{corollary}
		Every triangle-free planar graph is $1$-defective   online $3$-DP-colorable.
	\end{corollary}
 
We denote a $(2,1)$-decomposition of a graph $G$   by a pair $(D,M)$, where $M$ is a matching in $G$, and $D$ is an acyclic orientation of $G-M$ with maximum out-degree $\Delta^+(D) \le 2$.

		All the graphs considered in this paper are simple and finite. For a graph $G$  and a vertex $x \in V(G)$, $N_G(u)$ is the set of neighbors of $u$, $d_G(u) = |N_G(u)|$ is the degree of $u$, $\delta(G), \Delta(G)$ are the minimum and maximum degree of vertices in $G$, respectively. For $S \subseteq  V(G)$, $N_G(S)= \cup_{v \in S}N_G(v)$, and $G-S$ is the graph obtained from $G$ by deleting vertices in $S$ and all the edges incident with them. For a set $E$ of unordered pairs of vertices,
		$G+E$ (resp. $G-E$)  is the graph obtained from $G$ by adding (resp. deleting) the elements of $E$ to the edge set (resp. from the edge set) of $G$. If  $E=\{uv\}$, then we write  $G+uv$ (resp. $G-uv$) for $G+E$ (resp. $G-E$).
		For a digraph $D$ and a set $A$ of ordered pairs on $V(D)$, define $D+A$, $D-A$, $D+(u,v)$, $D-(u,v)$ similarly. For two vertices $x,y \in V(D)$, let $D-xy=D-\{(x,y),(y,x)\}$.  For two graphs (or digraphs) $G$ and $H$, the union $G \cup H$ of $G$ and $H$ has vertex set  $ V(G) \cup V(H)$ and edge set $E(G) \cup E(H)$. We denote by $d_D^+(v)$ and $d_D^-(v)$ the out-degree and the in-degree of $v$ in $D$, respectively. Note that a graph is $d$-degenerate if and only if $G$ has an acyclic orientation $D$ with maximum out-degree at most $d$. 
  
  		Assume $H$ is a proper subgraph of a connected graph $G$. For a component  with vertex-set $K$ of $G-V(H)$, the subgraph of $G$ induced by    $K \cup N_G(K)$    is called a \emph{bridge of $H$ in $G$}. For a bridge $B$ of $H$ in $G$, the elements of $V(H) \cap V(B)$  are called its \emph{vertices of attachment to $H$.} An edge not in $H$ linking two vertices of $H$ is a \emph{singular bridge of $H$ in $G$}, and is also called \emph{chord of $H$} when $H$ is a cycle. Similarly, a $2$-chord of a cycle $C$ in $G$ is a bridge with three vertices and two vertices of attachment to $C$.

  	If $G$ is a plane graph and $C$ is  a 
  	cycle in $G$, then ${\rm Int}(C)$ is the subgraph of $G$ induced by all vertices   inside or on $C$, ${\rm Ext}(C)$ is the subgraph of $G$ induced by  all vertices   outside or on $C$.   We denote by $B_G$  the boundary walk of the infinite face of $G$.    For two vertices $u$ and $v$ in  $B_G$,  $B_G[u,v]$ is the path on $B_G$ from $u$ to $v$ in clockwise direction.  Vertices and edges in $B_G$ are called  \emph{boundary vertices} and  \emph{boundary edges} of $G$. If $uv$ is a boundary edge, then $v$ is called a \emph{boundary neighbor} of $u$. If $uv$ is a chord of $B_G$, then $v$ is called a \emph{chord neighbor} of $u$. For a subset $S$ of $B_G$,  denote by $T_G(S)$ the chord neighbors of vertices in $S$ in $G$. If $G$ is $2$-connected, then $B_G$ is a cycle, in this case, for each vertex $v \in B_G$,  we let $v^-, v^+$ be the previous and next boundary vertex of $v$ in $B_G$ in the clockwise direction respectively.
 
	\begin{definition}
		A {\em configuration} is a pair $(G, P)$, where $G$ is a 
		connected triangle-free plane graph and $P=v_1v_2v_3v_4$ is a path consisting of four consecutive vertices of  the boundary walk  of  $G$. 
	\end{definition}

	\begin{definition}
		\label{def-(2,1)-deco}
 	Assume $(G, P)$ is a configuration, where $P=v_1v_2v_3v_4$. Let $H$ be the block of $G$ containing edge $v_2v_3$. 
		Assume $\vec{a} =  a_1 a_2 a_3 a_4 \in \{0,1,2 \}^4$ and $\vec{b} = b_1 b_2 b_3 b_4 \in \{0,1 \}^4$. An $(\vec{a}, \vec{b})$-decomposition  of $(G,P)$ is a  $(2,1)$-decomposition  $(D,M)$ of $G-v_2v_3$ such that
		\begin{enumerate}[label= {(\arabic*)}]  
			\item \label{C1} For $i =1,2,3,4$,  $d_M(v_i) \le b_i$, and $d_D^+(v_i) \le a_i$.
			\item \label{C2} For any $v \in B_G \setminus \{v_1,v_2,v_3,v_4\}$, $d_D^+(v) \le 1$.
		\end{enumerate} 
		An $(\vec{a}, \vec{b})$-decomposition $(D,M)$ of $(G,P)$ is {\em relaxed}    if     (2) is replaced with 
		
		(2)': For any $v \in B_G\setminus T_H(\{v_2,v_3\})$, $d_D^+(v) \le 1$.

	 \noindent 
	\end{definition}

 	The difference between ``relaxed" and ``non-relaxed" version is that in the former, the chord neighbors of $v_2$ or $v_3$ in $H$ are allowed to have out-degree $2$ in $D$.   
  
	We write
	$$(G, v_1 v_2 v_3 v_4) \in \mathcal{G}(\vec{a}, \vec{b}), \ \text{ or }  \ (G, v_1 v_2 v_3 v_4) \in \mathcal{G}^*(\vec{a}, \vec{b})$$ if $(G,v_1 v_2 v_3 v_4)$ has an $(\vec{a}, \vec{b})$-decomposition, or $(G,v_1 v_2 v_3 v_4)$ has a relaxed $(\vec{a}, \vec{b})$-decomposition, respectively.

	

	
	Theorem \ref{main-thm0} follows from the following result.
	
	\begin{theorem}
		\label{thm-main1}
	 	For any configuration $(G,P)$,  
		$$(G, wxyz) \in \mathcal{G}(1001,1001).$$
	\end{theorem}

		\section{Special configurations}
		\label{configuration}
		
Theorem \ref{thm-main1} is proved by induction on the number of vertices. For the purpose of using induction, we need to prove a more technical statement. Roughly speaking, the more technical statement asserts that if $(G,P)$ does not contain certain configurations, then we can further require that in the decomposition $(D,M)$ of $(G,P)$, the two end vertices of $P$ are not covered by edges in $M$.

Now  we define some special configurations. 
For the configurations $(G, P)$ below, the path $P$ will be denoted by $wxyz$, and we assume that the four vertices $wxyz$ are in clockwise cyclic order, i.e., $w^+=x, x^+=y$ etc.. When there are more than one configurations involved, then the configurations will be denoted by $(G_i, w_ix_iy_iz_i)$.   For a statement $S$, we write $ (G, wxyz) \in \mathcal{G}(\vec{a}, \vec{b})|_{S}$ to mean that $(G, wxyz)$ has an $(\vec{a}, \vec{b})$-decomposition $(D,M)$ in which statement $S$ is true.   For example, $ (G, wxyz) \in \mathcal{G}(\vec{a}, \vec{b})|_{zz^+ \in M}$   means that $(G, wxyz)$ has an $(\vec{a}, \vec{b})$-decomposition $(D,M)$ with $zz^+ \in M$.   
	
	\begin{definition}
		\label{def-operation}
		Assume   $(G_1,w_1x_1y_1z_1), (G_2, w_2x_2y_2z_2)$ are two configurations. 
		\begin{itemize}
			\item  Denote by $$(G_1,w_1x_1y_1z_1) \oplus (G_2, w_2x_2y_2z_2) $$ the configurations $(G, w_1x_1y_1z_2)$,    where $G$ is obtained from the disjoint union of $G_1, G_2$ by identifying $x_2$ with $z_1$, $y_2$ with $y_1$.
			\item Denote by $$(G_1,w_1x_1y_1z_1) \hat{\oplus} (G_2, w_2x_2y_2z_2)$$ the configurations $(G, w_1x_1uz_2)$,    where $G$ is obtained from  $(G_1,w_1x_1y_1z_1) \oplus (G_2, w_2x_2y_2z_2) $ by adding vertex $u$ and edges $ux_1,uz_2$.
			\item Denote by $$(G_1,w_1x_1y_1z_1) \tilde{\oplus} (G_2, w_2x_2y_2z_2)$$  the configurations $(G, x_1uvz_2)$,  where $G$ is obtained from  $(G_1,w_1x_1y_1z_1) \oplus (G_2, w_2x_2y_2z_2)$ by adding vertices $u,v$ and edges $uv, ux_1,vz_2$.
		\end{itemize}
	\end{definition}
	
\begin{figure}[h]
		\centering
		\begin{minipage}[t]{0.33\textwidth}
			\centering
			\begin{tikzpicture}[>=latex,
			roundnode/.style={circle, draw=black,fill=black, very thick, minimum size=1mm, inner sep=0pt}]  
			\node [roundnode] (w1) at (0,1.5){};
			\node [roundnode] (x1) at (1,1.5){};
			\node [roundnode] (y1) at (2,1.5){};	
			\node [roundnode] (z1) at (2,0){};
			
			\node [roundnode] (w2) at (3.2,0){};
			\node [roundnode] (x2) at (2.2,0){};
			\node [roundnode] (y2) at (2.2,1.5){};	
			\node [roundnode] (z2) at (3.2,1.5){};
			
			\node at (0.2,1.25){$w_1$};
			\node at (1,1.25){$x_1$};
			\node at (1.8,1.25){$y_1$};
			\node at (1.9, -0.2){$z_1$};
			\node at (1, 0.75){$G_1$};
			
			\node at (3.5,0){$w_2$};
			\node at (2.4,-0.25){$x_2$};
			\node at (2.4,1.25){$y_2$};
			\node at (3.4, 1.25){$z_2$};
			\node at (3, 0.75){$G_2$};
			
			\draw (w1)--(x1)--(y1)--(z1);
			\draw (w2)--(x2)--(y2)--(z2);
			\end{tikzpicture}
			\subcaption{$(G_1, P_1) \oplus (G_2, P_2)$}
		\end{minipage}
		\begin{minipage}[t]{0.33\textwidth}
			\centering
			\begin{tikzpicture}[>=latex,	
			roundnode/.style={circle, draw=black,fill=black, very thick, minimum size=1mm, inner sep=0pt}]  
			\node [roundnode] (w1) at (0,1.5){};
			\node [roundnode] (x1) at (1,1.5){};
			\node [roundnode] (y1) at (2,1.5){};	
			\node [roundnode] (z1) at (2,0){};
			
			\node [roundnode] (w2) at (3.2,0){};
			\node [roundnode] (x2) at (2.2,0){};
			\node [roundnode] (y2) at (2.2,1.5){};	
			\node [roundnode] (z2) at (3.2,1.5){};
			\node [roundnode] (u) at (2,2.3){};
			
			\node at (0.2,1.25){$w_1$};
			\node at (1,1.25){$x_1$};
			\node at (1.8,1.25){$y_1$};
			\node at (1.9, -0.2){$z_1$};
			\node at (1, 0.75){$G_1$};
			\node at (2, 2.05){$u$};
			
			\node at (3.5,0){$w_2$};
			\node at (2.4,-0.25){$x_2$};
			\node at (2.4,1.25){$y_2$};
			\node at (3.4, 1.25){$z_2$};
			\node at (3, 0.75){$G_2$};
			
			\draw (w1)--(x1)--(y1)--(z1);
			\draw (w2)--(x2)--(y2)--(z2);
			\draw (x1)--(u)--(z2);
			\end{tikzpicture}
			\subcaption{$(G_1, P_1) \hat{\oplus} (G_2, P_2)$}
		\end{minipage}
		\begin{minipage}[t]{0.3\textwidth}
			\centering
			\begin{tikzpicture}[>=latex,	
			roundnode/.style={circle, draw=black,fill=black, very thick, minimum size=1mm, inner sep=0pt}]  
			\node [roundnode] (w1) at (0,1.5){};
			\node [roundnode] (x1) at (1,1.5){};
			\node [roundnode] (y1) at (2,1.5){};	
			\node [roundnode] (z1) at (2,0){};
			
			\node [roundnode] (u) at (1.5,2.3){};
			\node [roundnode] (v) at (2.5,2.3){};
			
			\node [roundnode] (w2) at (3.2,0){};
			\node [roundnode] (x2) at (2.2,0){};
			\node [roundnode] (y2) at (2.2,1.5){};	
			\node [roundnode] (z2) at (3.2,1.5){};
			
			\node at (0.2,1.25){$w_1$};
			\node at (1,1.25){$x_1$};
			\node at (1.8,1.25){$y_1$};
			\node at (1.9, -0.2){$z_1$};
			\node at (1, 0.75){$G_1$};
			
			\node at (1.6, 2.05){$u$};
			\node at (2.4, 2.05){$v$};
			
			\node at (3.5,0){$w_2$};
			\node at (2.4,-0.25){$x_2$};
			\node at (2.4,1.25){$y_2$};
			\node at (3.4, 1.25){$z_2$};
			\node at (3, 0.75){$G_2$};
			
			\draw (w1)--(x1)--(y1)--(z1);
			\draw (w2)--(x2)--(y2)--(z2);
			\draw (x1)--(u)--(v)--(z2);
			\end{tikzpicture}
			\subcaption{$(G_1, P_1) \tilde{\oplus} (G_2, P_2)$}
		\end{minipage}
		\caption{Illustration of the three operations.}
		\label{fig-operations}
	\end{figure}
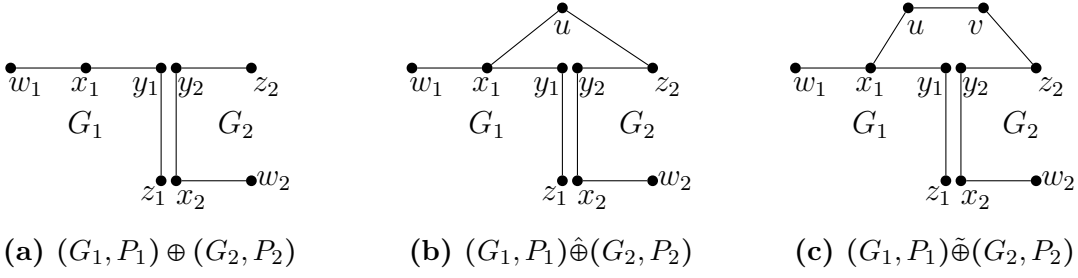

	Fig. \ref{fig-operations} illustrates the three operations. 
	 
	Now we define three families of special configurations $\mathcal{R}$, $\mathcal{Q}$ and $\mathcal{P}$.

	\begin{definition} \label{def-PR}
		The families of
		$\mathcal{R}$, $\mathcal{Q}$  and $\mathcal{P}$-configurations are defined recursively  as follows.  
		\begin{enumerate}[label= (\arabic*), ref=R\theenumi]
			\item \label{Rule1}  If $G $ is a $4$-cycle (respectively, a $5$-cycle)  and 
			$w, x, y, z$ are four consecutive boundary vertices of $G$,    
			then $(G,wxyz) \in \fr$ (respectively, $(G,wxyz)  \in \fp$).
			\item \label{Rule2} If $(G_1,w_1x_1y_1z_1) \in \fr, (G_2, w_2x_2y_2z_2) \in \fr$, then  
			$$ (G_1,w_1x_1y_1z_1) \oplus (G_2, w_2x_2y_2z_2) \in \fr,$$ 
			$$ (G_1,w_1x_1y_1z_1) \hat{\oplus} (G_2, w_2x_2y_2z_2) \in \fr,$$
			$$(G_1,w_1x_1y_1z_1) \tilde{\oplus} (G_2, w_2x_2y_2z_2) \in \fp.$$
			\item \label{Rule3}   If $(G_1,w_1x_1y_1z_1) \in \fp , (G_2, w_2x_2y_2z_2) \in \fr$, then  
			 $$(G_1,w_1x_1y_1z_1)  \hat{\oplus} (G_2, w_2x_2y_2z_2) \in \fp,$$ 
			 $$(G_1,w_1x_1y_1z_1)  \hat{\oplus} (G_2, w_2x_2y_2z_2) \in \fp.$$ 
		\end{enumerate} 
	\end{definition}
  
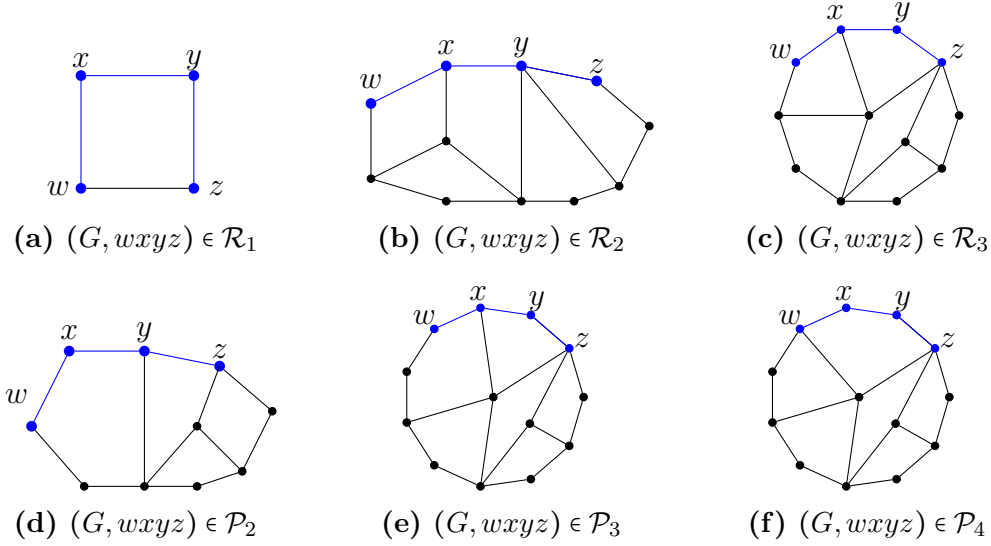
\begin{figure}
		\centering
		\begin{minipage}[t]{0.3\textwidth}
			\centering
			\begin{tikzpicture}[>=latex,	
			roundnode/.style={circle, draw=blue,fill=blue, very thick, minimum size=1mm, inner sep=0pt}]  
			\node [roundnode] (w) at (0,0){};
			\node [roundnode] (x) at (0,1.5){};
			\node [roundnode] (y) at (1.5,1.5){};	
			\node [roundnode] (z) at (1.5,0){};
			\node at (-0.3,0){$w$};
			\node at (0,1.7){$x$};
			\node at (1.5,1.7){$y$};
			\node at (1.8, 0){$z$};
			\draw[blue](w)--(x)--(y)--(z);
			\filldraw  (w)--(z);
			\end{tikzpicture}
			\subcaption{$(G, wxyz) \in \fr_1$}
		\end{minipage}
		\begin{minipage}[t]{0.3\textwidth}
		\centering
		\begin{tikzpicture}[>=latex,	
			roundnode/.style={circle, draw=blue,fill=blue, very thick, minimum size=1mm, inner sep=0pt}] 
			\node [roundnode] (w) at (0,1.3){};
			\node [roundnode] (x) at (1,1.8){};
			\node [roundnode] (z) at (3,1.6){};
			\node at (0,1.6){$w$};
			\node at (1,2.05){$x$};
			\node at (2,2.05){$y$};
			\node [roundnode] (y) at (2,1.8){};
			\node at (3, 1.8){$z$};
			\filldraw (0,0.3)circle(0.05)--(1,0)circle(0.05)--(2,0)circle(0.05)--(2.7,0)circle(0.05)--(3.3,0.2)circle(0.05)--(3.7,1)circle(0.05)--(z)--(y);
			\draw[blue] (w)--(x)--(y)--(z);
			\filldraw  (w)--(0,0.3);
			\draw (2,0)--(y)--(3.3,0.2);	
			\filldraw (2,0)--(1,0.8)circle(0.05)--(x);
			\draw (1,0.8)--(0,0.3);
		\end{tikzpicture}
		\subcaption{$(G, wxyz) \in \fr_2$}
	\end{minipage}
		\begin{minipage}[t]{0.3\textwidth}
			\centering
			\begin{tikzpicture}[>=latex,	
			bluenode/.style={circle, draw = blue, fill=blue, minimum size=1mm, inner sep=0pt},
			blacknode/.style={circle, draw = black, fill=black, minimum size=1mm, inner sep=0pt},
			] 
			\node [bluenode] (A1) at (36:1.2cm){};
			\node [bluenode] (A2) at (72:1.2cm){};
			\node [bluenode] (A3) at (108:1.2cm){};
			\node [bluenode] (A4) at (144:1.2cm){};
			\node [blacknode] (A5) at (180:1.2cm){};
			\node [blacknode] (A6) at (216:1.2cm){};	
			\node [blacknode] (A7) at (252:1.2cm){};
			\node [blacknode] (A8) at (288:1.2cm){};
			\node [blacknode] (A9) at (324:1.2cm){};
			\node [blacknode] (A10) at (360:1.2cm){};
			\node [blacknode] (u) at  (324:0.6cm){}; 
			\node  at (36:1.45cm){$z$};
			\node  at (72:1.45cm){$y$};
			\node  at (108:1.45cm){$x$}; 
			\node  at (144:1.45cm){$w$};
			\node [blacknode] (w) at (0:0cm){};
			\foreach \t in {1,3,5,7}{
				\draw (w) -- (A\t);}
			\draw (A7)--(u)--(A1);
			\draw (A9)--(u);
			
			\draw [blue] (A1)--(A2)--(A3)--(A4);
			\draw (A4)--(A5)--(A6)--(A7)--(A8)--(A9)--(A10)--(A1);
			\end{tikzpicture}
			\subcaption{$(G, wxyz) \in \fr_{3}$}
		\end{minipage}

		\medskip
		
		\begin{minipage}[t]{0.3\textwidth}
			\centering
			\begin{tikzpicture}[>=latex,	
				roundnode/.style={circle, draw=blue,fill=blue, very thick, minimum size=1mm, inner sep=0pt}] 
				\node [roundnode] (xn) at (1,1.8){};
				\node [roundnode] (z) at (3,1.6){};
				\node [roundnode] (w) at (0.5,0.8){};
				\node at (0.3,1.2){$w$};
				\node at (1,2.05){$x$};
				\node at (2,2.05){$y$};
				\node [roundnode] (x1) at (2,1.8){};
				\node at (3, 1.8){$z$};
				\filldraw   (w)--(1.2,0)circle(0.05)--(2,0)circle(0.05)--(2.7,0)circle(0.05)--(3.3,0.2)circle(0.05)--(3.7,1)circle(0.05)--(z);
				\draw[blue](z)--(x1);
				\draw[blue] (x1)--(xn);
			 
				\draw[blue] (w)--(xn);
				\draw (2,0)--(x1);	
				\draw (3.3,0.2)--(2.7,0.8)--(2,0);
				\filldraw (2.7,0.8)circle(0.05)--(3,1.6);
			\end{tikzpicture}
			\subcaption{$(G,wxyz) \in \fp_2$}
		\end{minipage}
		\begin{minipage}[t]{0.3\textwidth}
			\centering
			\begin{tikzpicture}[>=latex,	
			bluenode/.style={circle, draw = blue, fill=blue, minimum size=1mm, inner sep=0pt},
			blacknode/.style={circle, draw = black, fill=black, minimum size=1mm, inner sep=0pt},
			] 
			\node [bluenode] (A1) at (360/11:1.2cm){};
			\node [bluenode] (A2) at (720/11:1.2cm){};
			\node [bluenode] (A3) at (1080/11:1.2cm){};
			\node [bluenode] (A4) at (1440/11:1.2cm){};
			\node [blacknode] (A5) at (1800/11:1.2cm){};
			\node [blacknode] (A6) at (2160/11:1.2cm){};	
			\node [blacknode] (A7) at (2520/11:1.2cm){};
			\node [blacknode] (A8) at (2880/11:1.2cm){};
			\node [blacknode] (A9) at (3240/11:1.2cm){};
			\node [blacknode] (A10) at (3600/11:1.2cm){};
			\node [blacknode] (A11) at (0:1.2cm){};
			\node [blacknode] (u) at  (324:0.6cm){}; 
			\node  at (360/11:1.4cm){$z$};
			\node  at (720/11:1.4cm){$y$};
			\node  at (1080/11:1.4cm){$x$};
			\node  at (1440/11:1.4cm){$w$}; 
			\node [blacknode] (w) at (0:0cm){};
			\foreach \t in {1,3,6,8}{
				\draw (w) -- (A\t);}
			\draw (A8)--(u)--(A1);
			\draw (A10)--(u);
			\draw (A1)--(A2);
			\draw [blue] (A1)--(A2)--(A3)--(A4);
			\draw (A4)--(A5)--(A6)--(A7)--(A8)--(A9)--(A10)--(A11)--(A1);
			\end{tikzpicture}
			\subcaption{$(G, wxyz) \in \fp_{3}$}
		\end{minipage}
	\begin{minipage}[t]{0.3\textwidth}
		\centering
		\begin{tikzpicture}[>=latex,	
			bluenode/.style={circle, draw = blue, fill=blue, minimum size=1mm, inner sep=0pt},
			blacknode/.style={circle, draw = black, fill=black, minimum size=1mm, inner sep=0pt},
			] 
			\node [bluenode] (A1) at (360/11:1.2cm){};
			\node [bluenode] (A2) at (720/11:1.2cm){};
			\node [bluenode] (A3) at (1080/11:1.2cm){};
			\node [bluenode] (A4) at (1440/11:1.2cm){};
			\node [blacknode] (A5) at (1800/11:1.2cm){};
			\node [blacknode] (A6) at (2160/11:1.2cm){};	
			\node [blacknode] (A7) at (2520/11:1.2cm){};
			\node [blacknode] (A8) at (2880/11:1.2cm){};
			\node [blacknode] (A9) at (3240/11:1.2cm){};
			\node [blacknode] (A10) at (3600/11:1.2cm){};
			\node [blacknode] (A11) at (0:1.2cm){};
			\node [blacknode] (u) at  (324:0.6cm){}; 
			\node  at (360/11:1.4cm){$z$};
			\node  at (720/11:1.4cm){$y$};
			\node  at (1080/11:1.4cm){$x$}; 
			\node  at (1440/11:1.4cm){$w$};
			
			\node [blacknode] (w) at (0:0cm){};
			\foreach \t in {1,4,6,8}{
				\draw (w) -- (A\t);}
			\draw (A8)--(u)--(A1);
			\draw (A10)--(u);
			\draw (A1)--(A2);
			\draw [blue] (A1)--(A2)--(A3)--(A4);
			\draw (A4)--(A5)--(A6)--(A7)--(A8)--(A9)--(A10)--(A11)--(A1);
		\end{tikzpicture}
		\subcaption{$(G, wxyz) \in \fp_4$}
	\end{minipage}
		\caption{Examples of $\fr,\fp$-configurations.}
		\label{fig-configurations}
	\end{figure}

  For convenience,
  let $\fr_1,\fp_1$  denote the families of $\fr$-, $\fp$-configurations obtained by the first rule, respectively; let $\fr_2, \fp_2$ denote the corresponding configurations by the operation $\oplus$; let $\fr_2, \fp_2$ denote the corresponding configurations by the operation $\hat{\oplus}$; let $\fp_4$ denote the corresponding configurations by the operation $\tilde{\oplus}$. See Fig.~\ref{fig-configurations} for illustrations of these configurations. Note that all the graphs of these spacial configuration are $2$-connected. 
  First we have the following observation from the definition. 
	 
 	\begin{observation}
 		\label{obs-RP and Q}
 		$(G,wxyz) \in  \fr_3$ if and only if $(G-y,wxuz) \in \fr_2$, where $xyzu$ is  a $4$-face of $G$, and  $(G,wxyz) \in  \fp_4$ if and only if $(G-\{x,y\},w^-wvz) \in \fr_2$, where $wxyzv$ is   a $5$-face of $G$.
 	\end{observation}
 	
 	The operation $\oplus$ is not commutative, and $(G, wxyz) \in \fr$ does not implies that $(G,zyxw) \in \fr$. Nevertheless,  from the definition, we  know that the operation $\oplus$ is associative, 
 	$((G_1,w_1x_1y_1z_1) \oplus (G_2,w_2x_2y_2z_2)) \oplus (G_3,w_3x_3y_3z_3) = (G_1,w_1x_1y_1z_1) \oplus ((G_2,w_2x_2y_2z_2) \oplus (G_3,w_3x_3y_3z_3))$.
 
 	\begin{observation}\label{obs-associative} 
 		The following hold:
 		\begin{enumerate}[label= {(\arabic*)}] 
 			\item \label{obs1} If $(G,wxyz) \in  \fr_2$, then $(G,wxyz) = (G_1,w_1x_1y_1z_1) \oplus (G_2,w_2x_2y_2z_2)$ for some $(G_1,w_1x_1y_1z_1) \in  \fr$,  $(G_2,w_2x_2y_2z_2) \in \fr$.
 			\item \label{obs2} If $(G,wxyz) \in  \fr_3$, then $(G,wxyz) = (G_1,w_1x_1y_1z_1) \hat{\oplus} (G_2,w_2x_2y_2z_2)$ for some $(G_1,w_1x_1y_1z_1) \in \fr$, $(G_2,w_2x_2y_2z_2) \in \fr$. \qed
 		\end{enumerate}  
 	\end{observation}

 	\begin{observation}
 		\label{obs-key}
 		If $(G,wxyz) \in  \fr_i$, then $(G,z^+zyx) \in \fr_i$ for each $i \in \{1,2,3\}$.  If  
 		$(G,wxyz) \in  \fr_1, \fp_1, \fp_4$, then $(G,zyxw) \in  \fr_1, \fp_1, \fp_4$, respectively. 
 	\end{observation}
 	\begin{proof} The first part holds obviously if $(G,wxyz) \in \fr_1$. We assume that $(G,wxyz) \in \fr_2 \cup \fr_3$.  The proof is by induction. 
 		
 	First assume that $(G,wxyz) \in \fr_2$. By Observation \ref{obs-associative} and the definition of $\fr$-configurations, we may assume that $(G, wxyz) = (G_1,w_1x_1y_1z_1) \oplus (G_2,w_2x_2y_2z_2)$ for some $(G_1,w_1x_1y_1z_1) \in \fr$, and $(G_2,w_2x_2y_2z_2) \in \fr$. By induction hypothesis, $(G_2,z^+_2z_2y_2x_2) \in \fr$ and $(G_1,z^+_1z_1y_1x_1) \in \fr$, hence $(G, z^+zyx) = (G_2,z^+_2z_2y_2x_2) \oplus (G_1,z^+_1z_1y_1x_1) \in \fr_2$.
 	
 	If $(G,wxyz) \in \fr_3$, then by Observation \ref{obs-RP and Q}, $(G-y,wxuz) \in \fr_2$, where $xyzu$ is  a $4$-face of $G$. By induction, $(G-y,z^+,zux) \in \fr_2$, which implies that $(G,z^+zyx) \in \fr_3$. 
 		
 	The second part follows from the definition and the fact $(G,wxyz) \in \fr_2$ implying  $(G,z^+zyx) \in \fr_2$ as well as Observation \ref{obs-RP and Q}.
 	\end{proof}

 	The following lemma will be used in the proof of the main theorem. The proof of this lemma  is not difficult, but involves tedious verifications and is left to the last section. 
 	
 	 	\begin{lemma}
 	 		\label{main-lemmaA}
 	 		Let $$\tilde{\mathcal{G}} = \mathcal{G}(1101, 0000) \cap \mathcal{G}(1011, 0000) \cap \mathcal{G}(1002, 0000) \cap \mathcal{G}(2001, 0000).$$ Then  
 	 		\begin{enumerate}[label= {(\roman*)}]  
 	 			\item \label{lemR}  
 	 		 $\fr\setminus \fr_2 \subseteq   \tilde{\mathcal{G}} 
 	 			\cap \mathcal{G}(1001,1001)|_{zz^+\in M}  \cap \mathcal{G}(1001,1100) \cap \mathcal{G}(1001, 0011)|_{yz \in M}$. If  $(G,wxyz) \in  \fr_3$, then  $(G,wxyz) \in \mathcal{G}(1001,0001)|_{zz^+\in M}$.   
 	 			\item \label{lemP}
 	 			$\fp \setminus \fp_2 \subseteq    \tilde{\mathcal{G}}  \cap  \mathcal{G}(1001,0001)|_{zz^+\in M}\cap  \mathcal{G}(1001,1000)|_{w^-w
 	 				\in M}$. If $(G,wxyz) \in \fp_4$, then $(G,w^-wxy) \in   \mathcal{G}(1001,0000)$.  
  				
  				\item \label{lemR2}  $\fr_2 \subseteq  \mathcal{G}(1011,0000)|_{(z,y) \in D}$.
 	 		\end{enumerate} 
 	 	\end{lemma}

 	Assume $G$ is a  connected triangle-free plane graph, and $wxyz$ is a path of  boundary walk of $G$. We say  $G$ \emph{contains  an $\fr(xyz)$-configuration}, written as $G \sqsupset \fr(xyz)$, if  $G$ has a
 	subgraph $H$ such that $(H,uxyz)$ is an $\fr$-configuration  and $B_{H} \subseteq  B_G$ (for convenience, we also use $B_G$ to denote the  set  of boundary vertices of   $G$). Note that   $u \in B_G$ is not necessarily   $w$.  
 	Recall that each special configuration is 2-connected. 
 	So if $G$ contains  an $\fr(xyz)$-configuration $(H,uxyz)$, and $H'$ is the block of $G$ containing $xy$,  then $(H,uxyz)$ is contained in $H'$.  
 	
	We write   $G \nsqsupset \fr(xyz)$ if $G$ does not   contain any  $\fr(xyz)$-configuration.  Notions $G \sqsupset \fp(xyz)$ and $G \nsqsupset \fp(xyz)$ are defined similarly.

	Instead of proving Theorem \ref{thm-main1} directly, we prove the following more technical result.


	\begin{theorem}
	\label{main-thm}
	Let $G$ be a triangle-free plane graph, $wxyz$  a path in the boundary of $G$, and $H$ the block of $G$ containing edge $xy$.  Then the following hold:
	\begin{enumerate}[label= {(C\arabic*)}] 
		\item \label{M0}  $(G, wxyz)$ has a $(1001,1001)$-decomposition  $(D,M)$ such that $N_{M}(w),N_M(z) \subseteq B_G$.
		\item \label{M1} If $x$ or $y$ is incident to a chord of $B_H$, then 
		$(G, wxyz) \in \mathcal{G}^*(1001,0000)$.
		\item \label{M2} If  $G\nsqsupset \fr(xyz)$, $G \nsqsupset \fr(yxw)$, $G \nsqsupset \fp(xyz)$ and $G \nsqsupset \fp(yxw)$, then 
		$(G, wxyz) \in \mathcal{G}(1001,0000)$. 
		\item \label{M3} If $G \nsqsupset \fr(xyz)$, then $(G, wxyz) \in \mathcal{G}(1001,1000)$.
	\end{enumerate}
\end{theorem}

		\section{Proof of Theorem \ref{main-thm}}

		Assume  Theorem \ref{main-thm} is not true and $G$ is a counterexample   with minimum number of vertices. 
		
		 Note that   if $(G', P')$ is a smaller configuration, then $(G',P')$ is not a counterexample to Theorem \ref{main-thm}. Hence  $(G',P')$ has all the decompositions as stated in the conclusion. For brevity, we say `` $(D,M)$ is a  desired decomposition of $(G',P')$''   to mean that $(D,M)$ is any of these decompositions of $(G',P')$. 
		 To be precise,  $(D,M)$ is a $(1001,1001)$-decomposition with $N_M(w), N_M(w) \subseteq B_G$ for \ref{M0},   a relaxed   $(1001,0000)$-decomposition  for \ref{M1}, a $(1001,0000)$-decomposition  for \ref{M2}, a $(1001,1000)$-decomposition  for \ref{M3}.

		Now we shall prove a sequence of properties of $(G, wxyz)$  that  eventually lead to a contradiction.
		
		\begin{claim}
			\label{claim-min-deg}
			 Every vertex from $V(G)\setminus V(B_G)$ is of degree at least $3$.
		\end{claim}	
		\begin{proof}
			 If $d_G(v) \le 2$, then by the minimality of $G$, $(G-v,P)$ has a desired decomposition $(D,M)$. Let $D'=D+\{(v,u): u \in N_G(v)\}$. Then $(D',M)$ is  the desired decomposition of $(G,P)$.
		\end{proof}
		
		\begin{claim}
			\label{claim-2-connected}
			$G$ is $2$-connected.
		\end{claim}
		\begin{proof}
			Assume to the contrary, $G$ is not $2$-connected. Let $H_1,H_2,\ldots, H_s$ be the bridges of $H$ in $G$, and each  {$H_i$ has exactly}  one vertex  $u_i$ of attachment.  
			
		 Note that   $H_i-u_i$ contains at most one of $w$ and $z$.  If $H_i-u_i$ contains $w$, then $u_i=x$ and let $v_i=w$; If $H_i-u_i$ contains $z$, then $u_i=y$ and let $v_i=z$.  
		Otherwise let $v_i$ be an arbitrary boundary neighbor of $u_i$ in $v_i$.
		 
		 For each $i=1, \ldots, s$, $(H_i, u_i^-u_iv_iv_i^+)$ has a  $(1001,1001)$-decomposition $(D_i,M_i)$ (here $u_i^-,v_i^+$ refer to the corresponding neighbours on the boundary of the associated bridges).

		 Let $(D_H,M_H)$ be a $(1001,1001)$-decomposition (or a relaxed $(1001,0000)$-decomposition) of  $(H,x^-xyy^+)$ ($x^-$ and $y^+$ refer to neighbors of $x,y$ on $B_H$), let $D= D_H \cup \bigcup_{i=1}^{s}D_i+\bigcup_{i=1}^{s}(v_i,u_i)$, and $M=M_H \cup \bigcup_{i=1}^{s}M_i$.
		Then $(D,M)$ is a  $(1001,1001)$-decomposition (or a relaxed $(1001,0000)$-decomposition) of $(G,wxyz)$. This proves \ref{M0} and \ref{M1}.    
		 
	 Now we prove \ref{M2} and \ref{M3}.  If both $w,z$ are in $H$, then let $(D_{H},M_H)$ be a desired decomposition of $(H, wxyz)$.  If neither $w$ nor $z$ is in $H$, then let $(D_{H},M_H)$ be a   $(1001,1001)$-decomposition of $(G, x^-xyy^+)$. If $H=xy$, i.e., $xy$ is a cut edge, then $D_H=M_H=\emptyset$. Let $D=D_H \cup \bigcup_{i=1}^{s}D_i+\bigcup_{i=1}^{s}(v_i,u_i)$, and $M=M_H \cup \bigcup_{i=1}^{s}M_i$. It is easy to check that $(D,M)$ is a desired decomposition of $G$. 
		 
	Assume exactly one of $w,z$ is in $H$. 	   
		Without loss of generality, we may assume that $w \in V(H)$ but $z \notin V(H)$. In this case, $G \nsqsupset \fr(xyz)$ and $G \nsqsupset \fp(xyz)$ as $z$ is in the bridge of $H$ and is not the vertex of attachment. If $G \nsqsupset \fr(yxw)$, then by \ref{M3}, $(G,y^+yxw) \in \mathcal{G}(1001,1000)$. Let $(D'_{H},M'_H)$ be a $(1001,1000)$-decomposition of $(G,y^+yxw)$, $D'= D'_H \cup \bigcup_{i=1}^{s}D_i+\bigcup_{i=1}^{s}(v_i,u_i)$, and $M=M'_H \cup \bigcup_{i=1}^{s}M_i$.
		 Then $(D'_H,M'_H)$ is a  $(1001,0000)$-decomposition of $(G,wxyz)$, \ref{M2} and \ref{M3} hold. 
		 
		 Assume that $G \sqsupset \fr(yxw)$.  There is nothing to prove for \ref{M2}. As $z=v_i$ for some $i$, the decomposition $(D,M)$ constructed in the proof of  \ref{M0} is a  $(1001,1000)$-decomposition of  $(G, wxyz)$. This proves \ref{M3}.
		\end{proof}

		\begin{claim}
			\label{claim-no-separating cycle}
			$G$ has no separating $4$- and $5$-cycles, and $B_G$ is neither a $4$-cycle nor a $5$-cycle.
		\end{claim}
		\begin{proof} 
			Suppose $C_0=v_1v_2\cdots v_p$ is a separating $4$- or $5$-cycle of $G$.  Let $G_1 = {\rm Ext}(C_0)$ and $G_2={\rm Int}(C_0)$. By the minimality of $G$, $(G_1,wxyz)$ has the desired    decomposition  $(D_0,M_0)$.
			We will extend $(D_0,M_0)$ to the whole graph $G$ as follows.
			
			First assume that $p=4$. We choose labels of $V(C_0)$ so that $v_1v_2 \notin M_0$, and without loss of generality, assume that $(v_2,v_1) \in D_0$. Let $G'=G_2- \{v_3\}$. By \ref{M0}, $(G', v^-_4v_4v_1v_2)$  has a   $(1001,1001)$-decomposition $(D',M')$ with $N_{M'}(v_2) \subseteq B_{G'}$. 
			Let 
			\[
			(D'',M'')=\begin{cases} (D' +  (u_1,v_2), M'- u_1v_2), &\text{ if $u_1v_2 \in M'$}, \cr 
			(D',M'), &\text{ otherwise.}
			\end{cases}
			\]
		 Let 
			 $D=D_0 \cup D'' +\{(u,v_3)|u \in N_{G_2}(v_3)\setminus \{v_2,v_4\}\}$ and $M=M_0 \cup M''$.  Note that if $u_1v_2 \in M'$, then $u_1$ is a boundary vertex of $G'$ and hence $d_{D'}(u_1) \le 1$, implying that $d_{D''}(u_1) \le 2$. Moreover, as $G$ is triangle-free, $u_1$ is not adjacent to $v_3$, and hence $d_{D}(u_1) \le 2$.    
			 Since $v_4,v_1$ are sinks in $D'$ and $v_2$ is a sink in $D'-(v_2,v_1)$, there is no directed path from $v_2$ to $v_1$ in $D$ even if $(v_1,v_4) \in D_0$. So $(D,M)$ is the desired decomposition. See the left of Fig.~\ref{fig:deco of claim3}  for illustrations.

			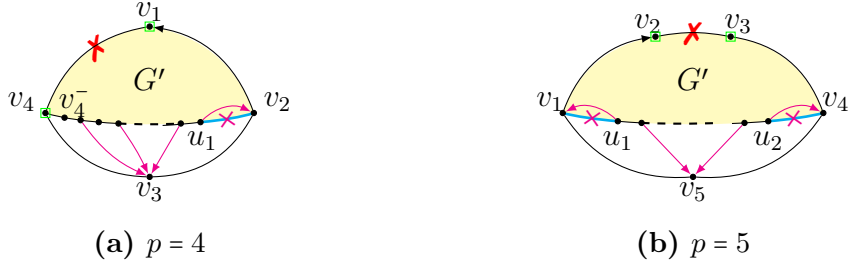
\begin{figure}[htbp]
				\centering
				\begin{minipage}[t]{0.45\textwidth} 
				\centering
				\begin{tikzpicture}[>=latex,	
					roundnode/.style={circle, draw=black,fill=black, minimum size=0.7mm, inner sep=0pt},
					greensquare/.style={rectangle, draw=green, minimum size=1.2mm, inner sep=0pt}]  
					
					\draw[draw=none, fill = yellow!30] ($(0,0)+(90:1)$) to [bend left] ($(0,0)+(330:1.6 and 0.3)$) arc (330:210:1.6 and 0.3) to [bend left] ($(0,0)+(90:1)$);
					
					\draw[cyan,line width=1pt] ($(0,0)+(330:1.6 and 0.3)$) arc (330:295:1.6 and 0.3); 
					
					\node [roundnode] (v1) at (90:1){};
					\node [roundnode] (v2) at ($(0,0)+(330:1.6 and 0.3)$){};	
					\node [roundnode] (v3) at (270:1){};
					\node [roundnode] (v4) at ($(0,0)+(210:1.6 and 0.3)$){};
					
					\node [roundnode] (u1) at ($(0,0)+(295:1.6 and 0.3)$){};
					\node [roundnode] (u2) at ($(0,0)+(285:1.6 and 0.3)$){};
					\node [roundnode] (u3) at ($(0,0)+(255:1.6 and 0.3)$){};
					\node [roundnode] (u4) at ($(0,0)+(245:1.6 and 0.3)$){};
					\node [roundnode] (u5) at ($(0,0)+(235:1.6 and 0.3)$){};
					\node [roundnode] (u6) at ($(0,0)+(225:1.6 and 0.3)$){};
					
					\draw (u1) arc (305:285:1.6 and 0.3); 	
					\draw (v4) arc (210:255:1.6 and 0.3); 
					\draw[dashed,line width =0.8pt] (u3) arc (255:285:1.6 and 0.3);
					
					\node at ($(-1,0)$){$v^-_4$}; 
					\node at ($(0,0)+(300:1.4 and 0.6)$){$u_1$}; 
					\node[rotate=10]  at ($(0,0)+(310:1.6 and 0.3)$){\textcolor{magenta}{\scriptsize \dmark}};

					\node at (90:1.2cm){$v_1$};
					\node at (0:1.7cm){$v_2$};
					\node at (270:1.2cm){$v_3$}; 
					\node at (180:1.7cm){$v_4$}; 
					
					\node at (0,0.3){$G'$}; 
					
					\node[rotate=30] at (135:1cm){\red{\xmark}};
					\draw[->] (v2) to [bend right]  (v1);
					\draw (v2) to [bend left] (v3);
					\draw (v3) to [bend left] (v4);
					
					\draw (v4) to [bend left] (v1); 
					\draw[magenta,->] (u1) to [bend left=35] (v2); 
					
					\draw[magenta,->] (u2) to (v3); 
					\draw[magenta,->] (u3) to  [bend left=5] (v3); 
					\draw[magenta,->] (u5) to [bend right=15] (v3); 
					
					\node [greensquare] at (v1){};
					\node [greensquare] at (v4){};
				\end{tikzpicture}
				\subcaption{$p=4$}
				\label{fig:claim-no sep-cycle-1}
			\end{minipage}
				\begin{minipage}[t]{0.45\textwidth} 
				\centering
				\begin{tikzpicture}[>=latex,	
					roundnode/.style={circle, draw=black,fill=black, minimum size=0.7mm, inner sep=0pt},
					greensquare/.style={rectangle, draw=green, minimum size=1.2mm, inner sep=0pt}]  
					
					\draw[draw=none, fill = yellow!30] ($(0,0)+(60:1)$) to [bend left] ($(0,0)+(330:2 and 0.3)$) arc (330:210:2 and 0.3) to [bend left] ($(0,0)+(120:1)$) to [bend left =10]($(0,0)+(60:1)$);
					
					\draw[cyan,line width=1pt] ($(0,0)+(330:2 and 0.3)$) arc (330:300:2 and 0.3); 
					\draw[cyan,line width=1pt] ($(0,0)+(210:2 and 0.3)$) arc (210:240:2 and 0.3); 
					
					\node [roundnode] (v3) at (60:1){};
					\node [roundnode] (v2) at (120:1){};
					\node [roundnode] (v4) at ($(0,0)+(330:2 and 0.3)$){};	
					\node [roundnode] (v5) at (270:1){};
					\node [roundnode] (v1) at ($(0,0)+(210:2 and 0.3)$){};
					
					\node [roundnode] (u1) at ($(0,0)+(300:2 and 0.3)$){};
					\node [roundnode] (u2) at ($(0,0)+(290:2 and 0.3)$){}; 
					\node [roundnode] (u3) at ($(0,0)+(250:2 and 0.3)$){};
					\node [roundnode] (u4) at ($(0,0)+(240:2 and 0.3)$){};
					
					\draw (u1) arc (300:285:2 and 0.3);
					\draw (u3) arc (250:240:2 and 0.3);   
					\draw[dashed,line width =0.8pt] (u3) arc (255:285:2 and 0.3);
					
					\node at ($(0,0)+(300:2 and 0.6)$){$u_2$}; 
					\node at ($(0,0)+(240:2 and 0.6)$){$u_1$}; 
					\node[rotate=10]  at ($(0,0)+(312:2 and 0.3)$){\textcolor{magenta}{\scriptsize \dmark}};
					\node[rotate=-10]  at ($(0,0)+(228:2 and 0.3)$){\textcolor{magenta}{\scriptsize \dmark}};

					\node at (60:1.2cm){$v_3$};
					\node at (120:1.2cm){$v_2$};
					\node at (0:1.9cm){$v_4$};
					\node at (270:1.2cm){$v_5$}; 
					\node at (180:1.9cm){$v_1$}; 
					
					\node at (0,0.3){$G'$}; 
					
					\node at (90:0.9){\red{\xmark}};
					\draw[->] (v1) to [bend left]  (v2);
					\draw (v2) to [bend left=10] (v3);
					\draw (v4) to [bend right] (v3);
					
					\draw (v4) to [bend left] (v5); 
					\draw (v5) to [bend left] (v1); 
					\draw[magenta,->] (u1) to [bend left=35] (v4); 
					\draw[magenta,->] (u4) to [bend right=35] (v1);
					
					\draw[magenta,->] (u2) to (v5); 
					\draw[magenta,->] (u3) to  (v5);  
					
					\node [greensquare] at (v2){};
					\node [greensquare] at (v3){};
				\end{tikzpicture}
				\subcaption{$p=5$}
				\label{fig:claim-no sep-cycle-2}
			\end{minipage}
				\caption{Illustrations for the proof of Claim \ref{claim-no-separating cycle}.  In all figures, a vertex with a green square box indicates it  is not covered by the matching, an edges with a red mark indicates it is the center two vertices of $P$ in $(G,P)$, hence will be deleted when giving the decomposition.}
				\label{fig:deco of claim3}
			\end{figure}
			
			Now assume that $p=5$. We  name the vertices of $V(C_0)$ so that $v_1v_2, v_3v_4 \notin M_0$. Without loss of generality, we may assume that $(v_1,v_2) \in D_0$. 
			For the edge $v_3v_4$, either $(v_3,v_4) \in D_0$ or $(v_4,v_3) \in D_0$. Let $G'=G_2 \setminus \{v_5\}$.  
			By \ref{M0}, $(G',v_1v_2v_3v_4)$ has a  $(1001,1001)$-decomposition $(D',M')$, with $N_{M'}(v_1), N_{M'}(v_4) \subseteq B_{G'}$. 
			Let 
			\[
			(D'',M'')=\begin{cases} (D' +  (u_1,v_1)+(u_2,v_4), M'- \{u_1v_1,u_2v_4\}), &\text{ if $u_1v_1, u_2v_4 \in M'$}, \cr 
			(D' +  (u_1,v_1), M'- u_1v_1), &\text{ if $u_1v_1 \in M'$, $d_{M'}(v_4)=0$}, \cr 
			(D' +  (u_2,v_4), M'- u_2v_4), &\text{ if $u_2v_4 \in M'$, $d_{M'}(v_1)=0$}, \cr 
			(D',M'), &\text{ if $d_{M'}(v_1)=d_{M'}(v_4)=0$.}
			\end{cases}
			\]
			 As $d_{D'}^+(v_3) = d_{M'}(v_3)=0$, we know that $(v_4,v_3) \in D''$.   Let $D= D_0 \cup (D''-(v_4,v_3) ) + \{(u,v_5)|u \in N_{G_2}(v_5) \setminus \{v_1,v_4\}\}$ and $M=M_0 \cup  M'' $. 
			 It is clear that $(D,M)$ is a desired  decomposition.  As $(v_1,v_2) \in D'' \cap D_0$ and $d^+_{D''}(v_4)=0$,  $d^+_D(v_1)=d^+_{D_0}(v_1)$, $d^+_D(v_4)=d^+_{D_0}(v_4)$, and   if $(v_3,v_4) \in D_0$, then it is not contained in a directed cycle.

			Assume $B_G=wxyz$   is a $4$-cycle. We assume that $v_3=w$, $v_4=x$, $v_1=y$ and $v_2=z$.  As $G$ contains no triangle,   none of $x$ and $y$ is incident to a chord in $B_G$. Also we know that $B_G$ is chordless, $(B_G,wxyz) \in \fr$, so $G \sqsupset \fr(xyz)$, $G \sqsupset \fr(yxw)$. We only need to show that  $(G, wxyz) \in \mathcal{G}(1001,1001)$ with $wz \in M$. 
			Let  {$(D'',M'')$} be the decomposition of  $(G-w, x^-xyz)$ as in the proof of $G$ having no separating $4$-cycle,  where $x^-$ refers the vertex on $B_{G-w}$.  Let $D=D'' + (w,x) + \{(u,w)|u \in N_{G}(w)\setminus \{x,z\}\}$ and $M = M'' + wz$. Then $(D,M)$ is the desired decomposition.

		 Assume $B_G=wxyza$ is a $5$-cycle. As $G \nsqsupset \fr(xyz)$, it suffices to show that $(G, wxyz) \in \mathcal{G}(1001,1000)|_{aw \in M}$. Let $(D'',M'')$ be the decomposition
   		as shown in the proof of the 5-separating cycle case. Let $D= D''- yz + (a,z)+ (z,y)+ \{(u,a)|u \in N_{G}(a)\setminus \{z,w\} \}$ and $M = M'' + aw$. Then  $(D,M)$ is the desired decomposition.
		\end{proof}

 		\begin{observation}\label{claim-contain-belong}
			Assume $H$ is a subgraph of $G$ and $B_H$ has no chord, $v_1v_2v_3v_4$ is a path in $B_H$. If $H \sqsupset \fr_i(v_2v_3v_4)$ or  $H \sqsupset \fp_j(v_2v_3v_4)$, then since $H$ has no chord and separating $4$- and $5$-cycles,    $(H,v_1v_2v_3v_4) \in \fr_i$ or $(H,v_1v_2v_3v_4) \in \fp_j$, respectively, for each $i \in\{1,3\}$, $j \in \{1,3,4\}$.  \qed
		\end{observation}

		\begin{claim}
			\label{claim-no-chord}
			$B_G$ has no chord.
		\end{claim}
		\begin{proof}
			Assume to the contrary, $uv$ is a chord of $B_G$, which divides $G$ into $G_1$ and $G_2$. Assume $G_1$   contains $w$.

			\medskip
			\noindent
			{\bfseries Case 1.} \emph{There exists a chord $uv$ with $\{u,v\} \cap \{x,y\} = \emptyset$.} 
			
			\medskip
			\noindent
			In this case, $N_{G}(\{x,y\})=N_{G_1}(\{x,y\})$.  
			Observe that $G_1$ contains $\fr(yxw)$-, $\fr(xyz)$-, $\fp(xyz)$-, and
			 $\fp(yxw)$-configurations if and only if $G$ does. Let $(D_1,M_1)$ be a desired decomposition  of $(G_1, wxyz)$.
			 By \ref{M0}, $(G_2, u^-uvv^+)$ has a $(1001,1001)$-decomposition $(D_2, M_2)$ (where $u^-,v^+$ are vertices on $B_{G_2}$).  
			 Then  $(D_1 \cup D_2, M_1 \cup M_2)$ is the the corresponding desired decomposition of $(G, wxyz)$.  
			
			\medskip
			\noindent
			{\bfseries Case 2.} \emph{There exists a chord $uv$ with $\{u,v\} \cap \{x,y\} = \{y\}$.}   
			
			\noindent 
			In this case, $y$ is incident to a chord of $G$. We may choose $v$ so that $G_2$ has no chord incident with $u$. The proof of \ref{M0} is the same as Case 1, we do not repeat and only focus on \ref{M1}-\ref{M3}.

			We first show \ref{M1}, i.e., $(G,wxyz) \in \mathcal{G}^*(1001,0000)$. By \ref{M0}, $(G_2, vyzz^+)$ has a $(1001,1001)$-decomposition $(D^1_2,M^1_2)$.
			
			 
			 If $x$ or $y$ has chord neighbors in $B_{G_1}$, then by \ref{M1},  $(G_1,wxyv) \in \mathcal{G}^*(1001,0000)$. Let $(D'_1,M'_1)$ be a relaxed $(1001,0000)$-decomposition of $(G_1,wxyv)$. Then $(D^1_1 \cup D^1_2, M^1_1 \cup M^1_2)$ is a relaxed $(1001,0000)$-decomposition of $(G,wxyz)$.

			Assume none of $x$ and $y$ has a chord neighbor in $B_{G_1}$. Then $B_{G_1}$ is chordless (by the assumption of Case 2).

			If $G_1 \nsqsupset \fr(xyv)$, $G_1 \nsqsupset \fr(yxw)$, $G_1 \nsqsupset \fp(xyv)$ and $G_1 \nsqsupset \fp(yxw)$, then  by \ref{M2}, $(G_1,wxyv)$ has a $(1001,0000)$-decomposition $(D^2_1,M^2_1)$. Otherwise, by Observation \ref{claim-contain-belong} and the fact that $B_{G_1}$ is chordless, we have that $(G,wxyz)$ or $(G,zyxw)$ is  in $\fr_1 \cup \fr_3 \cup \fp_1 \cup \fp_3 \cup \fp_4$. It  follows from Lemma~\ref{main-lemmaA} \ref{lemR},\ref{lemP} that $(G_1,wxyv)$ has a $(1002,0000)$-decomposition  $(D^2_1,M^2_1)$. In any case,  $(D^2_1 \cup D^1_2, M^2_1 \cup M^1_2)$ is a relaxed $(1001,0000)$-decomposition of $(G,wxyz)$.
			
			Now we show \ref{M2} and \ref{M3} together by considering whether $G_2\sqsupset \fr(vyz)$ or not. 
			
			If $G_2\nsqsupset \fr(vyz)$, then by \ref{M3}, $(G_2,v^-vyz)$ has a $(1001,1000)$-decomposition $(D^2_2,M^2_2)$. For the proof of \ref{M2}, we have the assumptions that $G \nsqsupset \fr(yxw)$ and $G \nsqsupset \fp(yxw)$, which implies that  $G_1 \nsqsupset \fr(yxw)$ and $G_1 \nsqsupset \fp(yxw)$. By \ref{M3}, $(G_1,vyxw)$ has a $(1001,1000)$-decomposition $(D^3_1,M^3_1)$.  Therefore, $(D^3_1 \cup D^2_2, M^3_1 \cup M^2_2)$ is a  $(1001,0000)$-decomposition of $(G,wxyz)$.
			For the proof of \ref{M3}, it is clear that $(D^4_1 \cup D^2_2, M^4_1 \cup M^2_2)$ is a  $(1001,1000)$-decomposition of $(G,wxyz)$, where $(D^4_1,M^4_1)$ is a $(1001,1001)$-decomposition of $(G_1,wxyv)$.

			Thus assume that $G_2\sqsupset \fr(vyz)$. 
			
			For the proof of \ref{M2}, we have that $G_1 \nsqsupset \fr(xyv)$  and $G_1 \nsqsupset \fp(xyv)$, for otherwise $G \nsqsupset \fr(xyz)$ or $G \nsqsupset \fp(xyz)$, a contradiction to the assumptions in \ref{M2}. Also by the assumptions that $G \nsqsupset \fr(yxw)$  and $G \nsqsupset \fp(yxw)$ implies that  that $G_1 \nsqsupset \fr(yxw)$ and $G_1 \nsqsupset \fp(yxw)$. Thus by \ref{M2}, $(G_1,wxyv)$ has a $(1001,0000)$-decomposition $(D^5_1,M^5_1)$. By \ref{M0}, $(G_2,vyzz^+)$ has a $(1001,1001)$-decomposition $(D^3_2,M^3_2)$, we know that $(D^5_1 \cup D^3_2, M^5_1 \cup M^3_2)$ is a $(1001,0000)$-decomposition of $(G,wxyz)$.
			
			For the proof of \ref{M3}, as before, the assumption implies that $G_1 \nsqsupset \fr(xyv)$, thus $(G_1,wxyv)$ has a $(1001,1000)$-decomposition $(D^6_1,M^6_1)$. Therefore,  $(D^6_1 \cup D^3_2, M^6_1 \cup M^3_2)$ is a $(1001,1000)$-decomposition of $(G,wxyz)$.

			\medskip
			\noindent
			{\bfseries Case 3.} \emph{For every chord $uv$, $\{u,v\} \cap \{x,y\} = \{x\}$.}
			
			Since the statements of \ref{M0}-\ref{M2} are symmetric for $w,x,y,z$, the proof  directly follows from Case 2. We only focus on \ref{M3}. The assumption implies that $G_2 \nsqsupset \fr(xyz)$, so $(G_2,vxyz)$ has a $(1001,1000)$-decomposition $(D^4_2,M^4_2)$. On the other hand, by \ref{M0}, $(G_1,wxvv^+)$ has a $(1001,1001)$-decomposition $(D^6_1,M^6_1)$. Therefore,  $(D^6_1 \cup D^4_2, M^6_1 \cup M^4_2)$ is a $(1001,1000)$-decomposition of $(G,wxyz)\sqsupset \fr(vyz)$. 
		\end{proof}

		\begin{claim}
			\label{claim-no-PR}
			$G$ contains none of \ $\fp(xyz)$-, $\fp(yxw)$-, $\fr(xyz)$- and \  $\fr(yxw)$-configurations.
		\end{claim} 
		\begin{proof}
		Assume to the contrary that $G$ contains $\fr(yxw)$-, $\fr(xyz)$-, $\fp(xyz)$- or $\fp(yxw)$-configurations. By Claim \ref{claim-no-chord}, $G$ contains none of $\fr_2(yxw)$-, $\fr_2(xyz)$-, $\fp_2(xyz)$- or $\fp_2(yxw)$-configurations. 
	By Observation~\ref{claim-contain-belong}, $(G,wxyz)$ or $(G,zyxw)$  is in $\fr_1 \cup \fr_3 \cup \fp_1 \cup \fp_3 \cup \fp_4$. In each case, by Lemma~\ref{main-lemmaA}\ref{lemR},\ref{lemP}, $(G,wxyz) \in \mathcal{G}(1001,1001)$ with $w^-w \in M$ if $d_{M}(w)=1$ and $zz^+ \in M$ if $d_{M}(z)=1$ (note that in the cases $(G,wxyz) \in \fr_1 \cup \fr_3$, $d_G(w)=2$, so if $d_{M}(w)=1$ and $d_{M}(x)=0$, then it must be that $w^-w \in M$). This proves \ref{M0}. For \ref{M3}, 
	if $(G, wxyz)$  or $(G, zyxw)$ is in $\fp_1 \cup \fp_3 \cup \fp_4$, then  by Lemma \ref{main-lemmaA}\ref{lemP}, $(G,wxyz) \in \mathcal{G}(1001,1000)|_{ww^- \in M}$. If $G \sqsupset \fr(yxw)$, as $G \nsqsupset \fr(xyz)$,  we have that $(G,zyxw) \in \fr_3$, and by the second part of Lemma \ref{main-lemmaA}\ref{lemR}, $(G,zyxw) \in \mathcal{G}(1001,0001)|_{ww^- \in M}$.   Thus we proved \ref{M3}.   There is nothing to prove for \ref{M1} and \ref{M2}, as the ``if'' parts are not satisfied. This completes   the proof of Claim \ref{claim-no-PR}.  
		\end{proof}

		By Claims \ref{claim-no-chord} and \ref{claim-no-PR},
		$B_G$ has no chord and 	$G$ contains none of \ $\fp(xyz)$-, $\fp(yxw)$-, $\fr(xyz)$- and \  $\fr(yxw)$-configurations. To prove Theorem \ref{main-thm}, 
		 we   need to show that $(G,wxyz) \in \mathcal{G}(1001,0000)$.

\begin{figure}[H]
			\centering
			\begin{minipage}[t]{0.48\textwidth} 
				\centering
				\begin{tikzpicture}[>=latex,	
				roundnode/.style={circle, draw=black,fill=black, minimum size=0.7mm, inner sep=0pt},
				greensquare/.style={rectangle, draw=green, minimum size=1.2mm, inner sep=0pt}]

				\draw ($(0,0.1)+(105:1.8 and 1.2)$) arc (105:75:1.8 and 1.2) to ($(0,0.1)+(75:1.8 and 1.2)$);
				
				\draw [magenta,<-] ($(0,0.1)+(105:1.8 and 1.2)$) arc (105:130:1.8 and 1.2) to  ($(0,0.1)+(130:1.8 and 1.2)$);
				
				\draw  [magenta,<-] ($(0,0.1)+(75:1.8 and 1.2)$) arc (75:50:1.8 and 1.2) to ($(0,0.1)+(50:1.8 and 1.2)$);
				
				\draw [fill=cyan!25] ($(-0.1,0)+(130:1.8 and 1.2)$) arc (130:240:1.8 and 1.2) to (-0.1,0.2) to ($(-0.1,0)+(130:1.8 and 1.2)$);
				\draw [fill=gray!25] ($(0.1,0)+(50:1.8 and 1.2)$) arc (50:-60:1.8 and 1.2) to (0.1,0.2) to ($(0.1,0)+(50:1.8 and 1.2)$);
				
				\draw [fill=yellow!30] ($(0,-0.1)+(240:1.8 and 1.2)$) arc (240:300:1.8 and 1.2) to (0,0.1) to ($(0,-0.1)+(240:1.8 and 1.2)$);	
				
				\node [roundnode] (w) at ($(0,0.1)+(130.5:1.8 and 1.2)$){}; 
				\node [roundnode] (x) at ($(0,0.1)+(105:1.8 and 1.2)$){};
				\node [roundnode] (y) at ($(0,0.1)+(75:1.8 and 1.2)$){};
				
				\node [roundnode] (z) at ($(0,0.1)+(50:1.8 and 1.2)$){};
				
				\node [roundnode] (x0) at ($(0.1,0)+(50:1.8 and 1.2)$){}; 
				\node [roundnode] (xi) at ($(0.1,0)+(300:1.8 and 1.2)$){};
				\node [roundnode] (xi-2) at ($(0,-0.1)+(300:1.8 and 1.2)$){};
				
				\node [roundnode] (xi+1) at ($(0,-0.1)+(240:1.8 and 1.2)$){};
				\node [roundnode] (xi+1') at ($(0,-0.1)+(260:1.8 and 1.2)$){};
				\node [roundnode] (xi+1-2) at ($(-0.1,0)+(240:1.8 and 1.2)$){};
				\node [roundnode] (br+1) at ($(-0.1,0)+(130:1.8 and 1.2)$){}; 
				
				\node [roundnode] (u-N) at (0,0.3){};
				\node [roundnode] (u-S) at (0,0.1){};
				\node [roundnode] (u-W) at (-0.1,0.2){};
				\node [roundnode] (u-E) at (0.1,0.2){};
				
				\node [rotate=-20]  at (-0.9,0.6){\red{\xmark}};
				
				\node [rotate=70]  at (-0.35,-0.5){\red{\xmark}}; 
				\node [rotate=10]  at (0.75,0.6){\red{\xmark}};

				\node  at (-1.2,0){$G_{L}$}; 
				\node  at (0.1,-0.8){$G_i$};
				\node  at (1.2,0){$G_{R}$}; 
				\node  at (0,0.5){$u$};  
				
				\node  at ($(0,0.1)+(130:1.8 and 1.5)$){$w$}; 
				\node  at ($(0,0.1)+(105:1.8 and 1.5)$){$x$};  
				\node  at ($(0,0.1)+(75:1.8 and 1.5)$){$y$}; 
				\node  at ($(0,0.1)+(50:1.8 and 1.5)$){$z$};	
				\node [rotate=-25] at ($(0.1,0)+(45:2.3 and 1.5)$){$(x_0)$};
				\node [rotate=15] at ($(0,-0.1)+(305:1.8 and 1.5)$){$x_i$};
				\node [rotate=-15] at ($(0,-0.1)+(235:1.8 and 1.5)$){$x_{i+1}$};
				\node [rotate=-5] at ($(0,-0.1)+(265:1.8 and 1.5)$){$x^-_{i+1}$};
				\node [rotate=35]  at ($(-0.1,0)+(143:2.2 and 1.6)$){$(x_{r+1})$};
				
				\draw[->] ($(-0.1,0)+(170:1.8 and 1.2)$) arc (170:140:1.8 and 1.2) to (br+1);
				\draw [->](xi+1-2)--(u-W);
				\draw [->](xi-2)--(u-S);
			 
				\draw [<-]($(0,-0.1)+(242:1.8 and 1.2)$) arc (242:260:1.8 and 1.2) to (xi+1');
				
				\draw [magenta,->](u-N)--(w);
				\draw [magenta,->](u-N)--(z);
				
				\node [greensquare] at (xi+1){};
				\node [greensquare] at (u-S){};
				\node [greensquare] at (u-E){};
				\node [greensquare] at (x0){};
				\node [greensquare] at (xi-2){};
				\end{tikzpicture}  
			\end{minipage}
			\begin{minipage}[t]{0.48\textwidth} 
				\centering
				\begin{tikzpicture}[>=latex,	
				roundnode/.style={circle, draw=black,fill=black, minimum size=0.7mm, inner sep=0pt},
				greensquare/.style={rectangle, draw=green, minimum size=1.2mm, inner sep=0pt}]

				\draw [fill=green!10] ($(0,0.1)+(130:1.8 and 1.2)$) arc (130:50:1.8 and 1.2) to (0,0.1) to ($(0,0.1)+(130:1.8 and 1.2)$);
				
				\draw [fill=cyan!25] ($(-0.1,0)+(130:1.8 and 1.2)$) arc (130:240:1.8 and 1.2) to (-0.1,0) to ($(-0.1,0)+(130:1.8 and 1.2)$);
				
				\draw [fill=gray!25] ($(0.1,0)+(50:1.8 and 1.2)$) arc (50:-60:1.8 and 1.2) to (0.1,0) to ($(0.1,0)+(50:1.8 and 1.2)$);
				
				\draw [fill=yellow!30] ($(0,-0.1)+(240:1.8 and 1.2)$) arc (240:299:1.8 and 1.2) to (0,-0.1) to ($(0,-0.1)+(240:1.8 and 1.2)$);	
				
				\node [roundnode] (w) at ($(0,0.1)+(100:1.8 and 1.2)$){}; 
				\node [roundnode] (x) at ($(0,0.1)+(75:1.8 and 1.2)$){};
				\node [roundnode] (y) at ($(0,0.1)+(50:1.8 and 1.2)$){};
				
				\node [roundnode] (xr+1) at ($(0, 0.1)+(130:1.8 and 1.2)$){};
				\node [roundnode] (xr-) at ($(-0.1, 0)+(160:1.8 and 1.2)$){};

				\node [roundnode] (x0) at ($(0.1,0)+(50:1.8 and 1.2)$){}; 
				\node [roundnode] (xi) at ($(0.1,0)+(300:1.8 and 1.2)$){};
				\node [roundnode] (z) at ($(0.1,0)+(20:1.8 and 1.2)$){};
				\node [roundnode] (xi-2) at ($(0,-0.1)+(300:1.8 and 1.2)$){};
				
				\node [roundnode] (xi') at ($(0,-0.1)+(280:1.8 and 1.2)$){}; 
				
				\node [roundnode] (xi+1) at ($(0,-0.1)+(240:1.8 and 1.2)$){};
				\node [roundnode] (xi+1-2) at ($(-0.1,0)+(240:1.8 and 1.2)$){};
				\node [roundnode] (xr+1-2) at ($(-0.1,0)+(130:1.8 and 1.2)$){}; 
				
				\node [roundnode] (u-N) at (0,0.1){};
				\node [roundnode] (u-S) at (0,-0.1){};
				\node [roundnode] (u-W) at (-0.1,0){};
				\node [roundnode] (u-E) at (0.1,0){};
				
				\node [rotate=-40]  at ($(0.1,0)+(35:1.8 and 1.2)$){\red{\xmark}};
				\node [rotate=-30]  at ($(0,0.1)+(65:1.8 and 1.2)$){\red{\xmark}};
				
				\node [rotate=-30]  at (-0.7,0.4){\red{\xmark}}; 
				\node [rotate=-60] at (0.38,-0.6){\red{\xmark}};

				\node  at (-1.2,0){$G_{L}$}; 
				\node  at (0.1,-0.8){$G_i$};
				\node  at (1.2,0){$G_{R}$}; 
				\node  at (0,0.8){$G_{r}$}; 
				\node  at (0,0.3){$u$};  
				
				\node  at ($(0,0.1)+(100:1.8 and 1.5)$){$w$};  
				\node  at ($(0,0.1)+(75:1.8 and 1.5)$){$x$}; 
				\node  at ($(0,0.1)+(50:1.8 and 1.5)$){$y$};	
				\node [rotate=-25] at ($(0.1,0)+(45:2.3 and 1.5)$){$(x_0)$};
				\node   at ($(0.1,0)+(15:2.4 and 1.8)$){$z(x^+_0)$};
				\node [rotate=15] at ($(0,-0.1)+(305:1.8 and 1.5)$){$x_i$};
				\node [rotate=-15] at ($(0,-0.1)+(235:1.8 and 1.5)$){$x_{i+1}$};
				\node [rotate=-5] at ($(0,-0.1)+(280:1.8 and 1.5)$){$x^+_{i}$};
				\node [rotate=35]  at ($(-0.1,0)+(130:2.2 and 1.6)$){$x_{r}$};
				\node [rotate=35]  at ($(-0.1,0)+(160:2.2 and 1.6)$){$x^-_{r}$};
				
				\node [greensquare] at (w){};
				\node [greensquare] at (x){};
				\node [greensquare] at (y){};
				\node [greensquare] at (u-N){};
				\node [greensquare] at (u-W){};
				\node [greensquare] at (z){};
				\node [greensquare] at (xi+1){};
				\node [greensquare] at (xr+1-2){};
				\node [greensquare] at (u-S){};
				\node [greensquare] at (x0){};
				\node [greensquare] at (xi-2){}; 
				
				\draw[->] ($(-0.1,0)+(170:1.8 and 1.2)$) arc (170:140:1.8 and 1.2) to (br+1);
				\draw [->](xi+1-2)--(u-W);
				\draw [->](xi+1)--(u-S);
				\draw [->](u-E)--(x0);
				\draw [<-]($(0,-0.1)+(298:1.8 and 1.2)$) arc (298:280:1.8 and 1.2) to (xi');
				
				\end{tikzpicture}  
			\end{minipage}
			\caption{Illustriation for Claim \ref{claim-no-2-chord-1} (left) and Claim \ref{claim-no-2-chord-2} (right)}
			\label{fig-Claim 5 and 6}
		\end{figure}
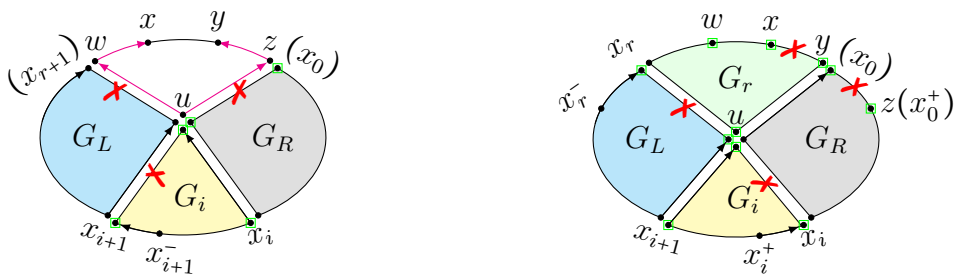

		\begin{claim}
			\label{claim-no-2-chord-1}
			$G$ has no $2$-chord of form $wuz$.
		\end{claim} 
		\begin{proof}
			Assume to the contrary there exist such a 2-chord $wuz$. Then $wuzyx$ is a facial cycle. Let $G' = G-\{x,y\}$. 
			
			 {First assume that  $B_{G'}$ has no chord incident to $u$. As $d_{G}(u) \ge 3$,  $G'$ is not a $4$-cycle. Also $G' \nsqsupset  \fr(wuz)$, for otherwise, by Observation \ref{claim-contain-belong}, $(G',w^-wuz) \in \fr$, and hence $(G,wxyz) \in \fp$, a contradiction to Claim \ref{claim-no-PR}.} By \ref{M3}, $(G',w^-wuz)$ has a $(1001,1000)$-decomposition $(D',M')$. Then  $(D'-(z,u)+(u,z)+(u,w)+(w,x)+(z,y), M')$ is a $(1001,0000)$-decomposition of $(G,wxyz)$. 
			
			Assume $G'$ has a chord incident to $u$. Let $x_0=z,x_1,x_2,\ldots x_r,x_{r+1}=w$ be the boundary neighbors of $u$ along the path $B_G-\{x,y\}$ in order. Let $C_i = B_G[x_i,x_{i+1}]\cup \{x_{i+1}ux_i\}$, $G_i = {\rm Int}(C_i)$. It is clear that $G_i$ has no chord incident to $x_{i+1}$ and $u$. 
			If $(G_i,x^-_{i+1}x_{i+1}ux_i) \in \fr$  for all $i=0,,1\ldots,r$, then  $(G,wxyz) \in \fp$, a contradiction. 
			
		So there exists an index $i$ such that  $G_i \nsqsupset \fr(x_{i+1},u,x_i)$. Let $C_L=B_G[x_{i+1},w]\cup \{wux_{i+1}\}$, $C_R=B_G[z,x_i]\cup \{x_iuz\}$, $G_L = {\rm Int}(C_L)$ and $G_R = {\rm Int}(C_R)$. Note that it is possible that $i=r$ or $i=0$, in which case $G_L$ or $G_R$ are empty graphs.

	 Since $G_i \nsqsupset \fr(x_{i+1}ux_i)$, by \ref{M3},    $(G_i,x^-_{i+1}x_{i+1}ux_i)$ has a $(1001,1000)$-decomposition $(D_i,M_i)$. By \ref{M0},	$(G_L,w^-wux_{i+1})$ has a $(1001,1001)$-decomposition  $(D_L,M_L)$,    and $ (G_R,x_iuzz^+) $ has a   $(1001,1001)$-decomposition  $(D_R,M_R)$.  Note that $(x_i,u)$ is an arc in both $D_i$ and $D_R$.  {Let $D=D_L \cup D_i \cup D_R + (u,w) + (w,x)+(u,z)+(z,y)$,}  $M=M_L \cup M_i \cup M_R$,  see Fig.~\ref{fig-Claim 5 and 6} (left) for illustration. It is clear that $(D,M)$ is a  $(1001,0000)$-decomposition of $(G,wxyz)$, a contradiction.
		\end{proof}
		
		\begin{claim}
			\label{claim-no-2-chord-2}
			None of $x$ and $y$ is incident to a $2$-chord of $B_G$
		\end{claim}
		\begin{proof}
			Assume the claim is not true. By symmetry, we may assume that $y$ is incident to a 2-chord $yuv$  of $B_G$. Let $x_0=y,x_1,x_2,\ldots, x_r$ be the boundary neighbors of $u$   in clockwise order. Let $C_i = B_G[x_i,x_{i+1}]\cup \{x_{i+1}ux_i\}$ (the indices are modulo $r$), $G_i = {\rm Int}(C_i)$. We may choose $u$ such that $G_r$ is as small as possible, so $G_r$ has no $2$-chord incident with $y$. 
			
			If $G_i \nsqsupset \fr(x_iux_{i+1})$ for some $i$, then let $C_L=B_G[x_{i+1},x_{r}]\cup \{x_{r}ux_{i+1}\}$, $C_R=B_G[y,x_i]\cup \{x_iuy\}$,  and let $G_L = {\rm Int}(C_L)$, $G_R = {\rm Int}(C_R)$. It is possible that $i=r-1$ or $i=0$, in which case $G_L$ or $G_R$ are empty graphs. By \ref{M3},  $(G_i,x^+_ix_iux_{i+1})$ has a $(1001,1000)$-decomposition, say $(D_i,M_i)$.  By \ref{M0},  $(G_L, x_{i+1}ux_{r}x^-_{r})$ and $(G_R,uyzz^+)$ have $(1001,1001)$-decompositions, say $(D_L,M_L)$ and $(D_R,M_R)$. 
			As $C_{r}$ is chordless,  by the same arguments as in Case 2 of the proof of Claim \ref{claim-no-chord},   $(G_{r},wxyu)$ has a $(1002,0000)$-decomposition, say $(D_{r},M_{r})$. Let $D=D_{r} \cup D_L \cup D_i \cup D_R + (z,y)$,  and $M=M_{r} \cup M_L \cup M_i \cup M_R$, see Fig.~\ref{fig-Claim 5 and 6} (right) for illustration. It is easy to check that $(D,M)$ is a  $(1001,0000)$-decomposition of $(G,wxyz)$, a contradiction.

			Assume $G_i \sqsupset \fr(x_iux_{i+1})$ for each $i \in \{0,1,\ldots,r-1\}$. As $B(G_i)$ is chordless, by Observation \ref{claim-contain-belong}, $(G_i,x^+_ix_iux_{i+1}) \in \fr$ for each $i$.
				
			If $|C_{r}|=4$, then $x_{r}=w$, hence $(G,zyxw) \in \fr$ by Definition \ref{def-PR}, a contradiction to Claim \ref{claim-no-PR}. If $|C_{r}|=5$, then $(G,yxww^-) \in \fp_4$. By the second half of Observation \ref{obs-key}, $(G,w^-wxy) \in \fp_4$.  Then by the second half of Lemma~\ref{main-lemmaA}\ref{lemP}, $(G,wxyz)$ has a $(1001,0000)$-decomposition, a contradiction.
				
				Assume $|C_{r}|\geq 6$. 
				Let $C_W=B_G[x_0,x_{r}]\cup \{x_{r}uy\}$, $G_W = {\rm Int}(C_W)$. As $(G_i,x^+_ix_iux_{i+1}) \in \fr$ for each $i \in \{0,1,\ldots,r-1\}$,   by Definition \ref{def-PR}, $(G_W,zyux_{r}) \in \fr_2$. By Lemma \ref{main-lemmaA}\ref{lemR2},  $(G_W,zyux_{r})$ has a $(1011,0000)$-decomposition $(D_W,M_W)$ with $(x_{r},u) \in D_{W}$, this implies that there is no directed path from $u$ to $x_{r}$ in $D_W$.

				Now we prove  that	$(G_{r},wxyu) \in \mathcal{G}(1001,0001)$.
				
				If $(G_{r}, wxyu) \in \fr_3 \cup \fp_3 \cup \fp_4$, then by   Lemma~\ref{main-lemmaA}\ref{lemR},\ref{lemP}, $(G_{r},wxyu) \in \mathcal{G}(1001,0001)$.
				 Assume $(G_{r}, wxyu) \notin \fr_3 \cup \fp_3 \cup \fp_4$. Recall that $B_{G_r}$ is chordless, so $G_{r}$ contains none of $\fr_2(xyu)$-,$\fr_2(yxw)$-, $\fp_2(xyu)$- and $\fp_2(yxw)$-configurations. If $G_r \sqsupset \fp_4(xyu)$ or $G_r \sqsupset \fp_4(yxw)$, then $(G_{r}, wxyu) \in  \fp_4$ by Observation \ref{claim-contain-belong} as $B_{G_r}$ is chordless. Thus $G_{r} \nsqsupset \fp_4(xyu)$ and $G_{r} \nsqsupset \fp_4(yxw)$. Recall that there is no 2-chord of form $yvw$, we know that $G_{r} \nsqsupset \fr_3(yxw)$ and $G \nsqsupset \fp_3(yxw)$. As $|C_{r}|\geq 6$, $G_r \nsqsupset \fr_1(yxw)$, $G_r \nsqsupset \fp_1(yxw)$, $G_{r} \nsqsupset \fr_1(xyu)$, $G_{r} \nsqsupset \fp_1(xyu)$. Therefore, $G_{r}$ contains none of $\fr(xyu)$-, $ \fr(yxw)$-, $\fp(xyu)$- and $\fp(yxw)$-configurations. By \ref{M2}, $(G_{r},wxyu) \in \mathcal{G}(1001,0000)$. In any case, we have $(G_{r},wxyu) \in \mathcal{G}(1001,0001)$.  
					  
				Let $(D_{r}, M_{r})$ be a $(1001,0001)$-decomposition of $(G_{r},wxyu)$. 
			As $d_{D_{r}}^+(y)=0$ and $d_{M_{r}}(y)=0$,   $(u,y) \in D_{r}$, and  either $ux_{r}\in M_{r}$ or $(x_{r},u) \in D_{r}$.

			If $ux_{r}\in M_{r}$, then let $D=D_{r}\cup (D_W-(x_{r},u))$ and $M=M_{r} \cup M_W$.
			As  $x_{r}$ is a sink in $D_W-(x_{r},u)$ and recall that there is no directed path from $u$ to $x_{r}$,  we know that $D$ is acyclic and $x_{r}$ has out-degree at most one. 
			 Otherwise $(x_{r},u)$ is an arc in both $D_{r}$ and $D_W$. Let $D=D_{r} \cup D_W$ and $M=M_{r} \cup M_W$. Again $D$ is acyclic and $x_{r}$ has out-degree at most one. In any case,  $(D,M)$ is a  $(1001,0000)$-decomposition of $(G, wxyz)$. 
			\end{proof}
			
			\medskip
			 
		We choose a sequence of vertices $w_1w_2\ldots$ as follows:
		$w_1=y$. Suppose $i \ge 1$ and the path $P_i=w_1w_2\ldots w_i$ is chosen, and $w_i$ is a boundary vertex. If $B_G$ has no 2-chord $w_iuv$ with $v \notin P_i$, then let $w_{i+1} = w_i^+$ and $P_{i+1}=w_1w_2\ldots w_iw_{i+1}$. Assume $B_G$ has a 2-chord $w_iuv$ with $v \notin P_i$. Such a  2-chord $w_iuv$ separates $G$ into two parts $G_1, G_2$ with $V(G_1) \cap V(G_2) = \{w_i,u,v\}$ and $P_i \subseteq V(G_1)$. Let $w_iuv$ be such a 2-chord so that $G_2$ is maximum. Let $w_{i+1}=u$, $w_{i+2}=v$ and  $P_{i+2}=w_1w_2\ldots w_iw_{i+1}w_{i+2}$.  
		Repeat this process, we shall eventually construct a cycle 
		$C^*=w_1w_2\cdots w_sw_1$ such that $w_1=y$, $w_s=x$.
		
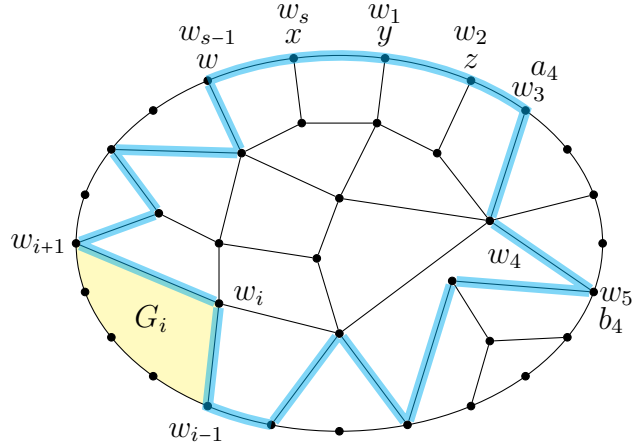
\begin{figure}[h]
			\centering  
			\begin{tikzpicture}[>=latex,
			roundnode/.style={circle, draw=black,fill=black, minimum size=1mm, inner sep=0pt}]
			
			\draw [draw=none, fill=yellow!30] ($(0,0)+(180:3.5 and 2.5)$) arc (180:240:3.5 and 2.5) to   (-1.6,-0.8) to ($(0,0)+(180:3.5 and 2.5)$);\draw ($(0,0)+(0:3.5 and 2.5)$) arc (0:360:3.5 and 2.5); 
			\node [roundnode] (w) at ($(0,0)+(120:3.5 and 2.5)$){}; 
			\node [roundnode] (x) at ($(0,0)+(100:3.5 and 2.5)$){};
			\node [roundnode] (y) at ($(0,0)+(80:3.5 and 2.5)$){};
			\node [roundnode] (z) at ($(0,0)+(60:3.5 and 2.5)$){};

			\foreach \i in {1,2,...,19}
			{\node [roundnode] (b\i) at ($(0,0)+(120+15*\i:3.5 and 2.5)$){};}
			
			\node [roundnode] (v1) at (-0.5,1.6){}; 
			\node [roundnode] (v2) at (0.5,1.6){}; 
			\node [roundnode] (v3) at (-1.3,1.2){}; 
			\node [roundnode] (v4) at (1.3,1.2){}; 
			\node [roundnode] (v5) at (-2.4,0.4){}; 
			\node [roundnode] (v7) at (-1.6,-0.8){}; 
			\node [roundnode] (v8) at (1.5,-0.5){};
			\node [roundnode] (v9) at (0,-1.2){};
			\node [roundnode] (v10) at (2,0.3){}; 
			\node [roundnode] (v11) at (2,-1.3){};

			\node [roundnode] (v13) at (-1.6,0){}; 
			\node [roundnode] (v15) at (-0.3,-0.2){};
			\node [roundnode] (v16) at (0,0.6){};

			\draw(x)--(v1)--(v2)--(y);
			\draw(w)--(v3)--(v1);
			\draw(z)--(v4)--(v2);
			\draw (v3)--(b2)--(v5)--(b4)--(v7)--(v9);
			\draw (v4)--(v10);
			\draw (v7)--(b8);
			\draw (v9)--(b9);
			\draw (v9)--(b11);
			
			\draw (v10)--(b19);
			\draw (b17)--(v10)--(b15);

			\draw (b15)--(v8)--(b11);
			\draw (v8)--(v11)--(b14); 
			\draw (v11)--(b12);
			\draw (v10)--(v16);
			
			\draw (v5)--(v13)--(v7);
			\draw (v9)--(v10);
			\draw (v16)--(v3)--(v13)--(v15)--(v9);
			\draw (v2)--(v16)--(v15);
			
			\draw[cyan,  draw opacity=0.5,line width = 4pt](w) arc (120:45:3.5 and 2.5)-- (v10)--(b15)--(v8)--(b11)--(v9)--(b9) arc(-105:-120:3.5 and 2.5)--(v7)--(b4)--(v5)--(b2)--(v3)--(w); 
			
			\node at ($(0,0)+(120:3.5 and 2.8)$){$w$}; 
			\node at ($(0,0)+(100:3.5 and 2.8)$){$x$}; 
			\node at ($(0,0)+(80:3.5 and 2.8)$){$y$}; 
			\node at ($(0,0)+(60:3.5 and 2.8)$){$z$}; 
			
			\node at ($(0,0)+(120:3.5 and 3.2)$){$w_{s-1}$}; 
			\node at ($(0,0)+(100:3.5 and 3.1)$){$w_s$}; 
			\node at ($(0,0)+(80:3.5 and 3.1)$){$w_1$}; 
			\node at ($(0,0)+(60:3.5 and 3.2)$){$w_2$};

			\node at ($(0,0)+(44:3.5 and 2.9)$){$w_3$}; 
			\node at ($(0,0)+(44:3.8 and 3.3)$){$a_4$}; 
			\node at ($(0,0)+(-14:3.8 and 2.9)$){$w_5$}; 
			\node at ($(0,0)+(-18:3.8 and 3.3)$){$b_4$};
			\node at (-2.5,-1){$G_i$}; 
			\node at (-1.2,-0.7){$w_{i}$}; 
			\node at (2.2,-0.2){$w_4$};

			\node at ($(0,0)+(240:3.8 and 2.9)$){$w_{i-1}$}; 
			\node at ($(0,0)+(180:4 and 2.9)$){$w_{i+1}$};

			\end{tikzpicture} 
			\caption{The blue cycle is $C^*$. Each vertex in $C^*$ has at least two neighbors in $B_G$ and each vertex in ${\rm Int}(C^*)$  has at most one neighbor in $B_G$. The yellow part is  $G_i$.}
			\label{fig:example of C*}
		\end{figure}

		Note that   each vertex $w_i$ either belongs to $B_G$ or $w_{i-1}w_iw_{i+1}$ is a 2-chord of $B_G$, and   each vertex in ${\rm Int}(C^*)\setminus V(C^*)$ is adjacent to  at most one vertex from $B_G$.
		See Fig.~\ref{fig:example of C*} for illustration. 	By Claim \ref{claim-no-2-chord-1} and Claim \ref{claim-no-2-chord-2}, $w_2=z$ and $w_{s-1}=w$.
		
		For   $w_i \in V(C^*)\setminus V(B_G)$, let $G_i={\rm Int}(B_G[w_{i-1},w_{i+1}]\cup \{w_{i+1}w_iw_{i-1}\})$.

		\begin{claim}
			\label{claim:C is not C*}
			$C^* \neq B_G$.
		\end{claim}
		\begin{proof}
			Assume to the contrary that $C^*=B_G$.   Let $G'= G-(V(C^*)\setminus\{x,y\})$.
			By \ref{M0}, $(G',x^-xyy^+)$ has a $(1001,1001)$-decomposition $(D',M')$ (here $x^-$ and $y^+$ are boundary vertices in $B_{G'}$). By the construction of $C^*$, each vertex $u \in B(G')-\{x,y\}$ has at most one neighbour in $B_G$.  By Claim \ref{claim-no-separating cycle}, $|B_G| \geq 6$, so there exists an edge $uv$ in $B_G$ such that $\{u,v\} \cap \{w,x,y,z\} =\emptyset$. Assume that $u,x,y,v$ appear in $B_G$ in this cyclic order. Let $D_P$ be the union of two directed path from $u$ to $x$ and $v$ to $y$. 
			Then let $M= M'+uv$ and $D = D_P \cup D'+\{(p,q)| p \in B(G')\setminus\{x,y\}, q \in V(B_G), pq \in E(G)\}$. Then  $(D,M)$ is a $(1001,0000)$-decomposition of $G$, a contradiction.		
		\end{proof}
 
\begin{figure}
			\centering
			\begin{minipage}[t]{0.45\textwidth} 
				\centering
				\begin{tikzpicture}[>=latex,	
				roundnode/.style={circle, draw=black,fill=black, minimum size=0.7mm, inner sep=0pt},
				greensquare/.style={rectangle, draw=green, minimum size=1.2mm, inner sep=0pt}]

				\draw [fill=cyan!25] ($(0,0.1)+(150:2.4 and 1.6)$) arc (150:30:2.4 and 1.6) to (0.7,0.1) to (-0.7,0.1) to ($(0,0.1)+(150:2.4 and 1.6)$);
				
				\draw [fill=gray!25] ($(-0.1,0)+(150:2.4 and 1.6)$) arc (150:240:2.4 and 1.6) to (-0.8,0) to ($(-0.1,0)+(150:2.4 and 1.6)$);
				\draw [fill=gray!25] ($(0.1,0)+(30:2.4 and 1.6)$) arc (30:-60:2.4 and 1.6) to (0.8,0) to ($(0.1,0)+(30:2.4 and 1.6)$);
				
				\draw [fill=yellow!30] ($(0,-0.1)+(240:2.4 and 1.6)$) arc (240:300:2.4 and 1.6) to (0.7,-0.1) to (-0.7,-0.1) to ($(0,-0.1)+(240:2.4 and 1.6)$);	
				
				\node [roundnode] (xn-1) at ($(0,0.1)+(127.5:2.4 and 1.6)$){}; 
				\node [roundnode] (xn) at ($(0,0.1)+(105:2.4 and 1.6)$){};
				\node [roundnode] (x1) at ($(0,0.1)+(75:2.4 and 1.6)$){};
				
				\node [roundnode] (x2) at ($(0,0.1)+(53.5:2.4 and 1.6)$){};
				
				\node [roundnode] (ai) at ($(0,0.1)+(30:2.4 and 1.6)$){};
				\node [roundnode] (ai2) at ($(0.1,0)+(30:2.4 and 1.6)$){};
				\node [roundnode] (bi1) at ($(0.1,0)+(300:2.4 and 1.6)$){};
				\node [roundnode] (bi2) at ($(0,-0.1)+(300:2.4 and 1.6)$){};
				\node [roundnode] (ai+1-1) at ($(0,-0.1)+(240:2.4 and 1.6)$){};
				\node [roundnode] (ai+1-2) at ($(-0.1,0)+(240:2.4 and 1.6)$){};
				\node [roundnode] (bi+1-1) at ($(-0.1,0)+(150:2.4 and 1.6)$){};
				\node [roundnode] (bi+1-2) at ($(0,0.1)+(150:2.4 and 1.6)$){};
				
				\node [roundnode] (wi+1-1) at (-0.8,0){};
				\node [roundnode] (wi+1-2) at (-0.7,0.1){};
				\node [roundnode] (wi+1-3) at (-0.7,-0.1){};
				\node [roundnode] (wi-1) at (0.8,0){};
				\node [roundnode] (wi-2) at (0.7,0.1){};
				\node [roundnode] (wi-3) at (0.7,-0.1){};
				
				\node [rotate=-20]  at (-1.5,0.33){\red{\xmark}};
				\node [rotate=0]  at ($(0,0)+(90:2.4 and 1.6)$){\red{\xmark}};
				\node  at (0,-0.154){\red{\xmark}}; 
				\node [rotate=20]  at (1.5,0.33){\red{\xmark}};
				\node [rotate=15]  at (1.05,-0.6){\scriptsize \textcolor{magenta}{\dmark}};	
				\node [rotate=-15]   at (-1.02,-0.6){\scriptsize \textcolor{magenta}{\dmark}};
				
				\node  at (-1.5,-0.2){$G_{j}$}; 
				\node  at (1.5,-0.2){$G_i$}; 
				\node  at (-0.4,0.35){$w_{j}$}; 
				\node  at (0.4,0.35){$w_{i}$}; 
				\node  at (0,-1){$G''$};
				\node  at (0,1.2){$G'$};
				
				\node  at ($(0,0.1)+(125:3 and 1.8)$){$w$}; 
				\node  at ($(0,0.1)+(102:3 and 1.8)$){$x$};  
				\node  at ($(0,0.1)+(78:3 and 1.8)$){$y$}; 
				\node  at ($(0,0.1)+(55:3 and 1.8)$){$z$};	
				\node  at ($(0,0.1)+(30:2.8 and 2)$){$w_{i-1}$};
				\node  at ($(0,-0.1)+(300:3 and 1.8)$){$w_{i+1}$};
				\node  at ($(0,-0.1)+(240:3 and 1.8)$){$w_{j-1}$};
				\node  at ($(-0.1,0)+(147:3 and 1.8)$){$w_{j+1}$};
				
				\node [greensquare] at (xn){};
				\node [greensquare] at (xn-1){};
				\node [greensquare] at (x1){};
				\node [greensquare] at (x2){};
				\node [greensquare] at (wi-1){};
				\node [greensquare] at (ai2){};	  
				\node [greensquare] at (bi+1-1){};
				\node [greensquare] at (wi+1-1){};
				\node [greensquare] at (wi+1-3){};
				\node [greensquare] at (wi-3){};
				\node [greensquare] at (bi2){};	
				\node [greensquare] at (ai+1-1){};	
				
			    \draw [->] (wi-2) to (0.8,0.6);
			    \draw [->] (0.8,0.6) to (0.3,0.8); 
			    \draw [dashed,->] (0.3,0.8) to (-0.3,0.8); 
			    \draw [->] (-0.3,0.8) to (-0.8,0.6); 
			    \draw [->] (-0.8,0.6) to (wi+1-2);  
				\draw [->](wi-3)--(0.3,-0.5);
				\draw [->](bi1)--(wi-1);
				\draw [->](ai+1-2)--(wi+1-1);
				\end{tikzpicture}  
			\end{minipage}
			\begin{minipage}[t]{0.45\textwidth} 
				\centering
				\begin{tikzpicture}[>=latex,	
				roundnode/.style={circle, draw=black,fill=black, minimum size=0.7mm, inner sep=0pt},
				greensquare/.style={rectangle, draw=green, minimum size=1.2mm, inner sep=0pt},
				bluesquare/.style={rectangle, draw=blue, minimum size=1.2mm, inner sep=0pt}]

				\draw [fill=cyan!25] ($(0,0.1)+(150:2.4 and 1.6)$) arc (150:30:2.4 and 1.6) to (0.6,0.1) to (-0.6,0.1) to ($(0,0.1)+(150:2.4 and 1.6)$);
				
				\draw [fill=gray!25] ($(-0.1,0)+(150:2.4 and 1.6)$) arc (150:220:2.4 and 1.6) to (-0.75,0.08) to ($(-0.1,0)+(150:2.4 and 1.6)$);
				\draw [fill=gray!25] ($(0.1,0)+(30:2.4 and 1.6)$) arc (30:-40:2.4 and 1.6) to (0.75,0.08) to ($(0.1,0)+(30:2.4 and 1.6)$);
				
				\draw [magenta,->] ($(0,-0.1)+(220:2.4 and 1.6)$) arc (220:240:2.4 and 1.6);
				\draw [magenta,->] ($(0,-0.1)+(300:2.4 and 1.6)$) arc (300:320:2.4 and 1.6);
				
				\draw [magenta,->] ($(0,-0.1)+(240:2.4 and 1.6)$) arc (240:259:2.4 and 1.6);
				\draw [magenta,->] ($(0,-0.1)+(260:2.4 and 1.6)$) arc (260:279:2.4 and 1.6);
				\draw [magenta,->] ($(0,-0.1)+(280:2.4 and 1.6)$) arc (280:299:2.4 and 1.6);
				
				\node [roundnode] (xn-1) at ($(0,0.1)+(127.5:2.4 and 1.6)$){}; 
				\node [roundnode] (xn) at ($(0,0.1)+(105:2.4 and 1.6)$){};
				\node [roundnode] (x1) at ($(0,0.1)+(75:2.4 and 1.6)$){};
				
				\node [roundnode] (x2) at ($(0,0.1)+(53.5:2.4 and 1.6)$){};
				
				\node [roundnode] (ai1) at ($(0,0.1)+(30:2.4 and 1.6)$){};
				\node [roundnode] (ai2) at ($(0.1,0)+(30:2.4 and 1.6)$){};
				\node [roundnode] (bi1) at ($(0.1,0)+(320:2.4 and 1.6)$){};
				\node [roundnode] (bi2) at ($(0,-0.1)+(320:2.4 and 1.6)$){};
				\node [roundnode] (aj1) at ($(0,-0.1)+(220:2.4 and 1.6)$){};
				\node [roundnode] (aj2) at ($(-0.1,0)+(220:2.4 and 1.6)$){};
				
				\node [roundnode] (bj1) at ($(-0.1,0)+(150:2.4 and 1.6)$){};
				\node [roundnode] (bj2) at ($(0,0.1)+(150:2.4 and 1.6)$){};
				
				\node [roundnode] (wj1) at (-0.6,0.1){}; 
				\node [roundnode] (wj2) at (-0.75,0.08){};
				\node [roundnode] (wi1) at (0.6,0.1){};
				\node [roundnode] (wi2) at (0.75,0.08){}; 
				
				\node [roundnode] (u1) at (0.4,0.1){};
				\node [roundnode] (u2) at (0.2,0.1){}; 
				\node [roundnode] (u3) at (0,0.1){};
				\node [roundnode] (u4) at (-0.2,0.1){}; 
				\node [roundnode] (u5) at (-0.4,0.1){}; 
				
				\node [roundnode] (v1) at ($(0,-0.1)+(300:2.4 and 1.6)$){};
				\node [roundnode] (v2) at ($(0,-0.1)+(280:2.4 and 1.6)$){};
				\node [roundnode] (v3) at ($(0,-0.1)+(260:2.4 and 1.6)$){};
				\node [roundnode] (v4) at ($(0,-0.1)+(240:2.4 and 1.6)$){}; 
				
				\draw [magenta,->](u1) to [bend left=15] (v2);
				\draw [magenta,->](u3) to [bend right=15] (v2);
				\draw [magenta,->](u4)--(v3);
				\draw [magenta,->](u5)--(v4);	
				
				\node [rotate=-30]  at (-1.5,0.35){\red{\xmark}};
				\node [rotate=0]  at ($(0,0)+(90:2.4 and 1.6)$){\red{\xmark}};
				
				\node [rotate=30]  at (1.5,0.35){\red{\xmark}};
				
				\node [rotate=-50]  at (1.3,-0.4){\scriptsize \textcolor{magenta}{\dmark}};	
				\node [rotate=50]  at (-1.3,-0.4){\scriptsize \textcolor{magenta}{\dmark}};
				
				\node  at (-1.8,-0.15){$G_{j}$}; 
				\node  at (1.8,-0.15){$G_i$}; 
				\node  at (-0.6,0.32){$w_{j}$}; 
				\node  at (0.6,0.32){$w_{i}$};
				
				\node  at (0,0.65){$G''$};
				
				\node  at ($(0,0.1)+(125:3 and 1.8)$){$w$}; 
				\node  at ($(0,0.1)+(102:3 and 1.8)$){$x$};  
				\node  at ($(0,0.1)+(78:3 and 1.8)$){$y$}; 
				\node  at ($(0,0.1)+(55:3 and 1.8)$){$z$};	
				\node  at ($(0,0.1)+(30:3 and 1.8)$){$w_{i-1}$};
				\node  at ($(0,-0.1)+(318:3 and 1.8)$){$w_{i+1}$};
				\node  at ($(0,-0.1)+(222:3 and 1.8)$){$w_{j-1}$};
				\node  at ($(-0.1,0)+(147:3 and 1.8)$){$w_{j+1}$}; 
				
				\node [greensquare] at (xn){};
				\node [greensquare] at (xn-1){};
				\node [greensquare] at (x1){};
				\node [greensquare] at (x2){};
				\node [greensquare] at (wi2){};
				\node [greensquare] at (ai2){};	 
				\node [greensquare] at (bj1){};
				\node [greensquare] at (wj2){}; 
				\draw [magenta,<-](aj1)--(wj1);
				\draw [magenta,->](wi1)--(bi2);
				\draw [->](aj2)--(wj2);
				\end{tikzpicture} 
				
			\end{minipage}
			\caption{Illustriation for Claim \ref{claim:no-consecutive wi} (left) and Claim \ref{claim:C*-bi=aj} (right)}
			\label{fig:claim 8 and 9} 
		\end{figure}
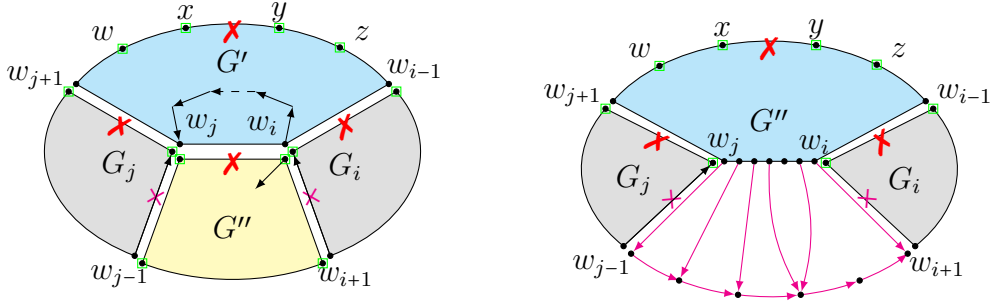

 In the rest of the proof, we shall frequently use the following observation.
 \begin{observation}
 	\label{obs-0000}
 	Assume $H$ is a  $2$-connected proper subgraph of $G$ with $|B_H| \geq 6$,  $v_1,v_2,v_3,v_4$ are four consecutive vertices in $B_H$. $B_H$ has no 2-chords of form $v_2vv_4$, $v_1vv_3$ and $v_1vv_4$.
 	If $T_{H}(v_2,v_3) = \emptyset$, then $(H,v_1v_2v_3v_4) \in \mathcal{G}(1001,0000)$. Otherwise,  $(H,v_1v_2v_3v_4) \in \mathcal{G}^*(1001,0000)$. 
 \end{observation}	
\begin{proof}
	If $B_H$ has a chord incident with $v_2$ or $v_3$, then by \ref{M1},  $(H,v_1v_2v_3v_4) \in \mathcal{G}^*(1001,0000)$. Assume $T_{H}(v_2,v_3) = \emptyset$.    
	 As $B_H$ has no 2-chord incident to $v_2$ or $v_3$ and no 2-chord joins $v_1$ and $v_4$, and  $|B_H| \geq 6$,
		 $H \nsqsupset  \fr(v_2v_3v_4)$, $H \nsqsupset  \fr(v_3v_2v_1)$, $H \nsqsupset  \fp(v_2v_3v_4)$, $H \nsqsupset   \fp(v_3v_2v_1)$. By   \ref{M2}, $(H,v_1v_2v_3v_4) \in \mathcal{G}(1001,0000)$.
\end{proof}
		\begin{claim}
			\label{claim:no-consecutive wi}
			$C^*$ has no chord.
		\end{claim}
		\begin{proof}
			Suppose that $w_i w_j$ is a chord of $C^*$ and $i < j$.
			By the construction of $C^*$, $w_i, w_j \notin V(B_G)$.  
			
			Let   $C'=B_G[w_{j+1},w_{i-1}]   \cup \{w_{i-1}w_iw_jw_{j+1}\}$ and $C''=B_G[w_{i+1}, w_{j-1}] \cup \{w_{j-1}w_{j}w_iw_{i+1}\}$, and let $G'={\rm Int}(C')$ and $G''={\rm Int}(C'')$.
			
			By the construction of $C^*$, $w_i$ and $w_j$ have no neighbors in $B_G[w_{i+1}^+,w^-_{j-1}]$.
			As $G$ is triangle-free,  neither $w_{i}w_{j-1}$ nor $w_{i+1}w_j$ is an edge of $G$. Hence $C''$ has no chord.

			Observe that by  Claim \ref{claim-no-2-chord-2}, neither $x$ nor $y$ is adjacent to $w_i$ or $w_j$, so $B_{G'}$ is chordless. Also, $B_{G'}$ has no $2$-chord   of form $wuz$, $wuy$, $xuz$ (by Claims \ref{claim-no-2-chord-1} and \ref{claim-no-2-chord-2}). By Observation \ref{obs-0000}, $(G', wxyz)$ has a  $(1001,0000)$-decomposition, say $(D',M')$. Without loss of generality, we may assume that either $(w_j,w_i) \in D'$ or $w_jw_i \in M'$.

			By \ref{M0}, we know that both $(G_i,w_{i+1}w_iw_{i-1}w^+_{i-1})$ and $(G_j,w_{j-1}w_jw_{j+1}w^-_{j+1})$ have  $(1001,1001)$-decompositions, say $(D_i,M_i)$ and $(D_j,M_j)$, respectively.   
			
		We claim that $(G'', w_{j-1}w_jw_iw_{i+1}) \in \mathcal{G}(1101,0000) \cap \mathcal{G}(1011,0000)$. Indeed,
		if  $G''$ contains none of the  $\fr(w_{j}w_iw_{i+1})$-, $\fr(w_iw_{j}w_{j-1})$-, $\fp(w_{j}w_iw_{i+1})$- and $\fp(w_iw_{j}w_{j-1})$-configurations, then by \ref{M2}, $(G'',w_{j-1}w_{j}w_iw_{i+1}) \in \mathcal{G}(1001,0000)$.  Otherwise, as $B_{G''}$ has  no chord,  by Observation \ref{claim-contain-belong},   $(G'',w_{j-1}w_{j}w_iw_{i+1})$ or $(G'',w_{i+1}w_iw_{j}w_{j-1})$ is in $\fr_1 \cup \fr_3 \cup \fp_1 \cup \fp_3 \cup \fp_4$. By Lemma \ref{main-lemmaA}\ref{lemR},\ref{lemP}, $(G'', w_{j-1}w_jw_iw_{i+1}) \in \mathcal{G}(1101,0000) \cap \mathcal{G}(1011,0000)$.
	 
		Let $(D_L,M_L)$ and  $(D_R,M_R)$ be a $(1101,0000)$-decomposition and a $(1011,0000)$-decomposition  of $(G'',w_{j-1}w_{j}w_iw_{i+1})$, respectively.
			
		 If there is a directed path from $w_i$ to $w_j$ in $D'$, then $(w_j, w_i) \notin D'$ and hence   $w_iw_j \in M'$. Let $M= M' \cup M_R \cup M_i \cup M_{j}$ and   $D= D'  \cup D_R \cup (D_i-(w_{i+1},w_i) \cup (D_{j}-(w_{j-1},w_j))$. Note that in $D_i-(w_{i+1},w_i)$, $w_{i+1}$ is a sink, and in $ D_{j}-(w_{j-1},w_j)$, $w_{j-1}$ is a sink. Hence $D$ is acyclic (see the left of Fig.~\ref{fig:claim 8 and 9} for illustration). The out-degree conditions are easily checked. So  $(D,M)$ is a  $(1001,0000)$-decomposition of $(G,wxyz)$.
			
		 If there is no directed path from $w_i$ to $w_j$ in $D'$, then let $M= M' \cup M_L \cup M_i \cup M_{j}$ and   $D= D'  \cup D_L \cup (D_i-(w_{i+1},w_i) \cup (D_{j}-(w_{j-1},w_j))$. 	
		 Again it is easy to check that $(D,M)$ is a  $(1001,0000)$-decomposition of $(G,wxyz)$, a contradiction. 	 
		\end{proof}
		
		\begin{claim}
			\label{claim:C*-bi=aj}
			If $w_i,w_j \in V(C^*)\setminus V(B_G)$, $i<j$, and $C^*[w_{i+1}, w_{j-1}] = B_G[w_{i+1}, w_{j-1}]$, then $j=i+2$ and  $w_{i+1}=w_{j-1}$.
		\end{claim}
		
		\begin{proof} 
			Assume to the contrary that $j > i+2$. Let  $C'=B_G[w_{j+1},w_{i-1}]\cup \{w_{i-1}w_iw_{i+1}\} \cup B_G[w_{i+1},w_{j-1}] \cup \{w_{j-1}w_jw_{j+1}\}$. Let $G'={\rm Int}(C')$   and $G''=G'-V(B_G[w_{i+1},w_{j-1}])$.  Clearly, $|B_{G''}| \geq 6$.
			
		 If $x$ or $y$ is incident to a chord of $  B_{G''}$, then by    \ref{M1},   $(G'', wxyz) \in \mathcal{G}^*(1001,0000)$.
		 
		 If none of $x$ and $y$ is incident to a chord of $B_{G''}$, then  as $|B_{G''}| \geq 6$ and $B_{G''}$ has no $2$-chord  of form $wuz$ (by Claim \ref{claim-no-2-chord-1}), $wuy$, $xuz$ (by Claim \ref{claim-no-2-chord-2}),  by Observation \ref{obs-0000},  $(G'', wxyz) \in \mathcal{G}(1001,0000)$.   Assume $(D'',M'')$ is   a relaxed or non-relaxed  $(1001,0000)$-decomposition of $(G'',wxyz)$ (depending on whether $x$ or $y$ is incident to a chord of $B_{G''}$ or not). Let  $(D_i,M_i)$ and $(D_j,M_j)$ be $(1001,1001)$-decompositions of $(G_i,w_{i+1}w_iw_{i-1}w^+_{i-1})$ and $(G_j,w_{j-1}w_jw_{j+1}w^-_{j+1})$, respectively.
			
		Orient the edges in $B_G[w_{i+1}, w_{j-1}]$ as a directed path  $D_P$  from $w_{j-1}$ to $w_{i+1}$  (see Fig.~\ref{fig:claim 8 and 9}). Let 
			$M = M''\cup M_i \cup M_j$, and 
			$D=D_P \cup D'' \cup (D_i-(w_{i+1},w_i)) \cup (D_j-(w_{j-1},w_j))+ \{(u,v)| \ u \in V(C'') \cap N_G(v), v \in V(B_G[w_{i+1},w_{j-1}]) \}$. 
			
			Note that if $x$ or $y$ has  a chord neighbor $v$ in $B_{G''}$, then since   $(D'',M'')$ is   a relaxed   $(1001,0000)$-decomposition of $(G'',wxyz)$, we may have  $d^+_{D''}(v)=2$. In this case,  
			by Claim \ref{claim-no-2-chord-2},    $v$  is not adjacent to any  vertex in $B_G[w_{i+1}, w_{j-1}]$. So $d^+_D(v)=d^+_{D''}(v)=2$. 
			
			Thus $(D,M)$ is a $(1001,0000)$-decomposition of $(G,wxyz)$, a contradiction.  
		\end{proof}

\begin{figure}[h]
		\centering
		\begin{minipage}[t]{0.45\textwidth} 
			\centering
				\begin{tikzpicture}[>=latex,	
				roundnode/.style={circle, draw=black,fill=black, minimum size=0.7mm, inner sep=0pt},
				greensquare/.style={rectangle, draw=green, minimum size=1.2mm, inner sep=0pt},
				bluesquare/.style={rectangle, draw=blue, minimum size=1.2mm, inner sep=0pt}]

				\draw [fill=cyan!25] ($(0,0.1)+(248:2.4 and 1.6)$) arc (248:80:2.4 and 1.6) to (0,0) to ($(0,0.1)+(248:2.4 and 1.6)$);

				\draw [fill=gray!25] ($(0,0)+(-70:2.4 and 1.6)$) arc (-70:-110:2.4 and 1.6) to (0,-0.2) to ($(0,0)+(-70:2.4 and 1.6)$);
				
				\draw [magenta,->] ($(0,0.1)+(-68:2.4 and 1.6)$) arc (-68:-41:2.4 and 1.6);
				\draw [magenta,->] ($(0,0.1)+(-40:2.4 and 1.6)$) arc (-40:-11:2.4 and 1.6);
				
				\draw [magenta,->] ($(0,0.1)+(-10:2.4 and 1.6)$) arc (-10:19:2.4 and 1.6);
				\draw [magenta,->] ($(0,0.1)+(20:2.4 and 1.6)$) arc (20:49:2.4 and 1.6);
				\draw [magenta,->] ($(0,0.1)+(50:2.4 and 1.6)$) arc (50:79:2.4 and 1.6);
				
				\node [roundnode] (xn-1) at ($(0,0.1)+(127.5:2.4 and 1.6)$){}; 
				\node [roundnode] (xn) at ($(0,0.1)+(105:2.4 and 1.6)$){};
				\node [roundnode] (x1) at ($(0,0.1)+(80:2.4 and 1.6)$){};
				  
				\node [roundnode] (x3) at ($(0,0.1)+(50:2.4 and 1.6)$){};
				\node [roundnode] (x4) at ($(0,0.1)+(20:2.4 and 1.6)$){};
				\node [roundnode] (x5) at ($(0,0.1)+(-10:2.4 and 1.6)$){};
				\node [roundnode] (x6) at ($(0,0.1)+(-40:2.4 and 1.6)$){};

				\node [roundnode] (wi1) at (0,0){};
				\node [roundnode] (wi2) at (0,-0.2){}; 
				
				\node [roundnode] (u1) at (0.07,0.27){};
				\node [roundnode] (u2) at (0.135,0.54){}; 
				\node [roundnode] (u3) at (0.2,0.81){};
				\node [roundnode] (u4) at (0.27,1.08){}; 
				\node [roundnode] (u5) at (0.34,1.35){}; 
				
				\node [roundnode] (w'i+1) at ($(0,0.1)+(248:2.4 and 1.6)$){};
				\node [roundnode] (w'i-1) at ($(0,0.1)+(-68:2.4 and 1.6)$){};
				\node [roundnode] (wi+1) at ($(0,0)+(-110:2.4 and 1.6)$){};
				\node [roundnode] (wi-1) at ($(0,0)+(-70:2.4 and 1.6)$){};
%
 				\draw [magenta,->](u4) to (x4);
				\draw [magenta,->](u2) to (x5);
				\draw [magenta,->](u1)--(x6);
				\draw [magenta,->](0,0)--(w'i-1);	
				
				\node [rotate=70]  at (-0.3,-0.8){\red{\xmark}};
				\node  at ($(0,0.1)+(90:2.4 and 1.6)$){\red{\xmark}};
			 	\node [rotate=-50]  at (0.35,-0.8){\scriptsize \textcolor{magenta}{\dmark}};

				\node  at (0,-1.2){$G_i$};   
				\node  at (-0.25,0.12){$w_i$};
				\node  at (0.1,1.35){$y'$};
				
				\node  at (-1,0){$G''$};
				
				\node  at ($(0,0.1)+(125:3 and 1.8)$){$w$}; 
				\node  at ($(0,0.1)+(102:3 and 1.8)$){$x$};  
				\node  at ($(0,0.1)+(78:3 and 1.8)$){$y$}; 
				\node  at ($(0,0.1)+(55:3 and 1.8)$){$z$};	 
				
				\node  at ($(0,0.1)+(249:3 and 2)$){$w_{i+1}$}; 
				\node  at ($(0,0.1)+(-71:3 and 2)$){$w_{i-1}$};

				\node [greensquare] at (xn){};
				\node [greensquare] at (xn-1){};
				\node [greensquare] at (wi+1){};
				\node [greensquare] at (0,-0.2){};
				\node [greensquare] at (x1){};
				\node [greensquare] at (0.34,1.35){}; 
				\draw [->](wi-1)--(0,-0.2);
			\end{tikzpicture}   
		\end{minipage}
		\begin{minipage}[t]{0.45\textwidth} 
			\centering
			\begin{tikzpicture}[>=latex,	
				roundnode/.style={circle, draw=black,fill=black, minimum size=0.7mm, inner sep=0pt},
				greensquare/.style={rectangle, draw=green, minimum size=1.2mm, inner sep=0pt},
				bluesquare/.style={rectangle, draw=blue, minimum size=1.2mm, inner sep=0pt}]

				\draw [fill=cyan!25] ($(0,0.1)+(208:2.4 and 1.6)$) arc (208:-28:2.4 and 1.6) to (0.75,0.1) to (-0.75,0.1) to ($(0,0.1)+(208:2.4 and 1.6)$);
				
				\draw [fill=gray!25] ($(-0.1,0)+(210:2.4 and 1.6)$) arc (210:270:2.4 and 1.6) to (-0.85,-0.05) to ($(-0.1,0)+(210:2.4 and 1.6)$);
				\draw [fill=gray!25] ($(0.1,0)+(-30:2.4 and 1.6)$) arc (-30:-90:2.4 and 1.6) to (0.85,-0.05) to ($(0.1,0)+(-30:2.4 and 1.6)$);

				\node [roundnode] (xn-1) at ($(0,0.1)+(60:2.4 and 1.6)$){}; 
				\node [roundnode] (xn) at ($(0,0.1)+(30:2.4 and 1.6)$){};
				\node [roundnode] (x1) at ($(0,0.1)+(0:2.4 and 1.6)$){};
				
				\node [roundnode] (x2) at ($(0,0.1)+(-28:2.4 and 1.6)$){};

				\node [roundnode] (wj1) at (-0.6,0.1){}; 
				\node [roundnode] (wj2) at (-0.75,0.08){};
				\node [roundnode] (wi1) at (0.6,0.1){};
				\node [roundnode] (wi2) at (0.75,0.08){}; 
				
				\node [roundnode] (u1) at (0.4,0.1){};
				\node [roundnode] (u2) at (0.2,0.1){}; 
				\node [roundnode] (u3) at (0,0.1){};
				\node [roundnode] (u4) at (-0.2,0.1){}; 
				\node [roundnode] (u5) at (-0.4,0.1){}; 
				
				\node [roundnode] (w6) at ($(-0.1,0)+(210:2.4 and 1.6)$){}; 
				\node [roundnode] (w2) at ($(0.1,0)+(-30:2.4 and 1.6)$){}; 
				\node [roundnode] (w4) at ($(0.1,0)+(-90:2.4 and 1.6)$){}; 
				\node [roundnode] (w4') at ($(-0.1,0)+(270:2.4 and 1.6)$){}; 
				
				\node [roundnode] (z+) at ($(0.1,0)+(-50:2.4 and 1.6)$){}; 
				\node [roundnode] (w6-) at ($(-0.1,0)+(230:2.4 and 1.6)$){}; 
				
				\node [roundnode] (w6') at ($(0,0.1)+(208.5:2.4 and 1.6)$){};

				\node [roundnode] (w3) at (0.85,-0.05){};
				\node [roundnode] (w5') at (-0.85,-0.05){};
				 
				\draw [magenta,->](u2)--(0.05,-1.5);
				\draw [magenta,->](u4)--(-0.05,-1.5);	
				
				\node [rotate=40]  at (-1.5,-0.5){\red{\xmark}};
				\node [rotate=-70]  at ($(0,0)+(15:2.4 and 1.6)$){\red{\xmark}};
				
				\node [rotate=-30]  at (1.5,-0.5){\red{\xmark}};
				
				\node [rotate=45]  at (-0.52,-0.8){\scriptsize \textcolor{magenta}{\dmark}};
				
				\node  at (-1,-0.8){$G_5$}; 
				\node  at (1,-0.8){$G_3$}; 
				\node  at (-0.6,0.32){$w_5$}; 
				\node  at (0.6,0.32){$w_3$};
				
				\node  at (0,0.8){$G''$};
				
				\node  at ($(0,0.1)+(60:3 and 1.8)$){$w$}; 
				\node  at ($(0,0.1)+(30:2.8 and 1.8)$){$x$};  
				\node  at ($(0,0.1)+(0:2.8 and 1.8)$){$y$}; 
				\node  at ($(0,0.1)+(-30:3.2 and 1.8)$){$z(w_2)$};
				\node  at ($(0,0.1)+(-50:3 and 2)$){$z^+$};  	
				\node  at ($(0,0.1)+(-90:3 and 2)$){$w_4$};  
				\node  at ($(0,0.1)+(210:2.8 and 2)$){$w_6$};
					\node  at ($(0,0.1)+(230:3 and 2)$){$w^-_6$};
				
				\node [greensquare] at (xn){};
				\node [greensquare] at (xn-1){};
				\node [greensquare] at (x1){};
				\node [greensquare] at (x2){}; 
				
				\node [greensquare] at (w3){};
				\node [greensquare] at (w6){};
				\node [greensquare] at (w2){};
				\node [greensquare] at (w5'){};
				\node [greensquare] at (w4){};
				
				\draw [magenta,->](wj2)--(w4');
			 
				\draw [->](w4')--(w5');
				\draw [->](wi2)--(x2);
			\end{tikzpicture}  
		\end{minipage}
		\caption{Illustriation for Claim \ref{claim-w3 and ws-2 is not in C} (left) and the end of the proof of the main theorem.}
		\label{fig:claim 11 and end}
	\end{figure}
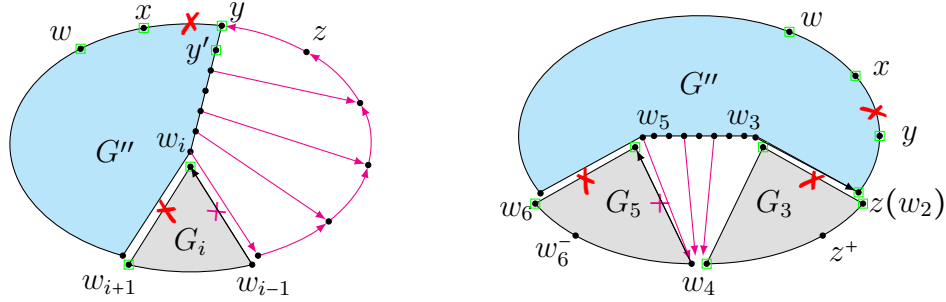

		\begin{claim}
			\label{claim-w3 and ws-2 is not in C}
			$w_3, w_{s-2} \notin V(B_G)$.
		\end{claim}
		\begin{proof} 
			Suppose to the contrary that $w_3 \in V(B_G)$. Let $i$ be the smallest index for which $w_i \notin V(B_G)$. By Claim \ref{claim:C is not C*}, such an index $i$ exists. Let $G'={\rm Int}(B_G[w_{i+1},w_{i-1}]\cup \{w_{i-1}w_iw_{i+1}\})$, and let $G''= G' \setminus V(B_G[z,w_{i-1}])$ (see the left of Fig. \ref{fig:claim 11 and end}). 
			
		 We first claim that $|B_{G''}| \geq 6$. For otherwise, by Claim \ref{claim-no-2-chord-2}, we have $|B_{G''}|= 5$, and $B_{G''}=wxyy'w_iw$  ($y'$ is the boundary neighbor of $y$ in $B_{G''}$ distinct to $x$). Also by Claim \ref{claim-no-2-chord-2}, $y'$ has no neighbor in $B_G[z,w_{i-1}]$, hence $d_{G}(y')=2$, a contradiction to Claim \ref{claim-min-deg}. 
			
		  By similar arguments as in Claim \ref{claim:C*-bi=aj}, $(G'',wxyy')$ has a relaxed   or non-relaxed $(1001,0000)$-decomposition   $(D'',M'')$. By definition of $G''$,  $y' \notin B_G$. 
			
			By \ref{M0}, $(G_i,w_{i-1}w_iw_{i+1}w^-_{i+1})$ has a   $(1001,1001)$-decomposition, say $(D_i,M_i)$. It follows that $(w_{i-1},w_i) \in D_i$. Let $D_P$ be the directed path $B_G[w_{i-1}, y]$ from $w_{i-1}$ to $y$. Let $M= M''\cup M_i$, and 
			$$D=D_P \cup D'' \cup (D_i-(w_{i-1},w_i)) + \{(u,v)| \ u \in V(B_{G''}) \cap N_G(v) , v \in V(B_G[w_{i-1},z]) \}.$$ 
			
			Then $(D,M)$ is a $(1001,0000)$-decomposition of $(G,wxyz)$, a contradiction. So $w_3 \notin V(B_G)$. 
			By symmetry, $w_{s-2} \notin V(B_G)$. 
		\end{proof}
	
		Recall that $w_2=z$. By Claim \ref{claim-w3 and ws-2 is not in C},  $w_4,w_5$ exist. By Claim \ref{claim-no-2-chord-1}, $w_4 \neq w$, 
		and by Claim \ref{claim:no-consecutive wi} and Claim \ref{claim:C*-bi=aj}, $w_4 \in V(B_G)$, $w_5 \notin V(B_G)$.
		 Let $G'={\rm Int}(B_G[w_6,z]\cup \{zw_3w_4w_5w_6\})$, and let $G''= G' \setminus \{w_4\}$.  By similar arguments as in Claim  \ref{claim:C*-bi=aj},
		  $(G'', wxyz)$ has a relaxed or non-relaxed $(1001,0000)$-decomposition $(D'',M'')$.  By Claim \ref{claim-no-2-chord-2},  
		  	any vertex $w \in B_{G''}[w_3,w_5]$ adjacent to $w_4$ 
		  	is not adjacent to $x$ or $y$ 
		  	(in particular, $w_3,w_5$ are not adjacent to $x$ or $y$).
		   Hence $d^+_{D''}(w_3)=1$ and $(w_3,z) \in D''$.  
		
		We claim that $(G_3,z^+zw_3w_4)\in \mathcal{G}(1011,1000)$.
		 If $G_3  \nsqsupset \fr(zw_3w_4)$, then by \ref{M3},
		  $(G_3,z^+zw_3w_4) \in \mathcal{G}(1001,1000)$. Otherwise, $G_3  \sqsupset \fr(zw_3w_4)$.  Note that all the possible chords in $B_{G_3}$ are incident with $w_3$, and $G_3$ has no separating $4$- and $5$-cycles, so it must be that $(G_3,z^+zw_3w_4) \in \fr$. By Lemma \ref{main-lemmaA}~\ref{lemR} and \ref{lemR2}, $(G,z^+zw_3w_4) \in \mathcal{G}(1011,0000)$.

		   Let $(D_3,M_3)$ be a $(1011,1000)$-decomposition of $(G_3, z^+zw_3w_4)$.
		 
		 By \ref{M0}, $(G_5,w^-_6w_6w_5w_4)$ has a $(1001,1001)$-decomposition $(D_5,M_5)$. Note  that $(w_4,w_5) \in D_5$.
		Let $M=M''\cup M_3 \cup M_5$, and  
		$$D = D'' \cup D_3\cup (D_5 -(w_4,w_5)) + \{(p,w_4)| \ p \in V(B_{G''})\cap N_{G}(w_4)  \}.$$
		Note that $w_4$ is not the tail of any arcs in the restriction to $G_5$ of $D$, hence $d^+_D(w_4) \leq 1$.
		Since $(w_3,z) \in D''$, $z$ is a sink in $D_3$, $(z,y) \in D'$ and $y$ is a sink in $D'$, there is no directed cycle in $D$.
	Then  $(D,M)$ is a $(1001,0000)$-decomposition of $(G,wxyz)$, a contradiction.
		
		This completes the proof Theorem \ref{main-thm}.

	\section{Proof of Lemma \ref{main-lemmaA}}
	
		\begin{definition}
			Assume $\vec{a}, \vec{b}, \vec{a}', \vec{b}',\vec{a}'', \vec{b}'' \in  \{0,1\}^4$. For $\circ \in \{ \oplus, \hat{\oplus}, \tilde{\oplus}\}$, we write 
			$$\mathcal{G}(\vec{a},\vec{b}) \circ \mathcal{G}(\vec{a}', \vec{b}') \subseteq  \mathcal{G}(\vec{a}'', \vec{b}''),$$
			if $(G_1, P_1) \in \mathcal{G}(\vec{a}, \vec{b})$ and 
			$(G_2, P_2) \in \mathcal{G}(\vec{a}', \vec{b})$ imply that $(G_1, P_1) \circ (G_2, P_2) \in \mathcal{G}(\vec{a}'', \vec{b}'')$.
		\end{definition}
		
		The following lemma lists some useful formulas.
		\begin{lemma}\label{lem-operation}
			$\mathcal{G}(\vec{a},\vec{b}) \circ \mathcal{G}(\vec{a}', \vec{b}') \subseteq  \mathcal{G}(\vec{a}'', \vec{b}'')$ in the following cases.
			\begin{enumerate}[label={(\roman*)}] 
				\item \label{formula-0}  $\mathcal{G}(1001,1001) \ \oplus \  \mathcal {G}(1001,1001)|_{zz^+ \in M} \subseteq   \mathcal{G}(1001,1001)|_{zz^+ \in M}$.
				\item \label{formula-1} For $\circ \in \{\oplus, \hat{\oplus},\tilde{\oplus}\}$, $\mathcal{G}(1011,0000) \ \circ \  \mathcal{G}(1001,1100) \subseteq   \mathcal{G}(1011,0000)$. 	
				\item \label{formula-2} $\mathcal{G}(1101,0000) \ \oplus \  \mathcal {G}(1001,1100) \subseteq   \mathcal{G}(1101,0000)$.  
				\item \label{formula-3} For $\circ \in \{\hat{\oplus},\tilde{\oplus}\}$, $\mathcal{G}(1001,0011) \ \circ \  \mathcal{G}(1001,1001)|_{zz^+ \in M} \subseteq   \mathcal{G}(1001,0001)|_{zz^+ \in M}$. 
				\item \label{formula-4} For $\circ \in \{\hat{\oplus},\tilde{\oplus}\}$, $\mathcal{G}(1011,0000) \ \circ \  \mathcal{G}(1001,1100) \subseteq   \mathcal{G}(1002,0000)$. 
				\item \label{formula-5} $\mathcal{G}(2001,0000) \  \hat{\oplus} \  \mathcal{G}(1001,1100) \subseteq   \mathcal{G}(2001,0000)$. 
				\item \label{formula-6}  $\mathcal{G}(1001,1000)|_{ww' \in M} \ \hat{\oplus} \  \mathcal {G}(1001,1100) \subseteq   \mathcal{G}(1001,1000)|_{ww' \in M}$. 
			\end{enumerate}
		\end{lemma}
		\begin{proof}
			Assume $(G_1, w_1x_1y_1z_1) \in \mathcal{G}(\vec{a}, \vec{b})$ and $(G_2, w_2x_2y_2z_2) \in \mathcal{G}(\vec{a}', \vec{b}')$. For each item,
			let $(D_1,M_1)$ be a corresponding $(\vec{a},\vec{b})$-decomposition of $(G_1, w_1x_1y_1z_1)$ and $(D_2, M_2)$ be a corresponding $(\vec{a}', \vec{b}')$-decomposition of $(G_2, w_2x_2y_2z_2)$.
			
			For \ref{formula-0} and \ref{formula-2}, let $D = D_1 \cup D_2$, and $M = M_1 \cup M_2$.
			
			For \ref{formula-1}, let $M = M_1 \cup M_2$. If $\circ = \oplus$, then let $D = D_1 \cup D_2$; If $\circ = \hat{\oplus}$, then let $D = D_1 \cup D_2 +(y_1,x_1)+(u,z_2)$; If $\circ = \tilde{\oplus}$, then let $D = D_1 \cup D_2 +(y_1,x_1)+(x_1,u)+(v,z_2)$.
			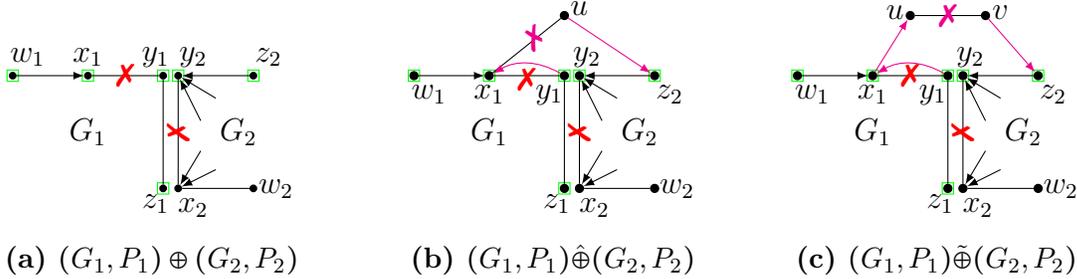
\begin{figure}[H]
				\centering
				\begin{minipage}[t]{0.33\textwidth}
					\centering
					\begin{tikzpicture}[>=latex,
					greensquare/.style={rectangle, draw=green, minimum size=1.5mm, inner sep=0pt},
					roundnode/.style={circle, draw=black,fill=black, very thick, minimum size=0.6mm, inner sep=0pt}]  
					\node [roundnode] (w1) at (0,1.5){};
					\node [roundnode] (x1) at (1,1.5){};
					\node [roundnode] (y1) at (2,1.5){};	
					\node [roundnode] (z1) at (2,0){};
					
					\node [roundnode] (w2) at (3.2,0){};
					\node [roundnode] (x2) at (2.2,0){};
					\node [roundnode] (y2) at (2.2,1.5){};	
					\node [roundnode] (z2) at (3.2,1.5){};
					
					\node [greensquare] (w1') at (0,1.5){};
					\node [greensquare] (x1') at (1,1.5){};
					\node [greensquare] (y1') at (2,1.5){};	
					\node [greensquare] (z1') at (2,0){};
					
					\node [greensquare] (y2') at (2.2,1.5){};	
					\node [greensquare] (z2') at (3.2,1.5){};
					
					\node at (0.2,1.75){$w_1$};
					\node at (1,1.75){$x_1$};
					\node at (1.9,1.75){$y_1$};
					\node at (1.9, -0.2){$z_1$};
					\node at (1, 0.75){$G_1$};
					
					\node at (3.5,0){$w_2$};
					\node at (2.4,-0.25){$x_2$};
					\node at (2.4,1.75){$y_2$};
					\node at (3.4, 1.75){$z_2$};
					\node at (3, 0.75){$G_2$};
					
					\draw [->] (w1)--(x1);
					\draw [->] (z2)--(y2);
					\draw [->] (2.5,0.5)--(x2);
					\draw [->] (2.7,0.25)--(x2);
					\draw [->] (2.5,0.9)--(y2);
					\draw [->] (2.7,1.2)--(y2);
					
					\draw (x1)--(y1)--(z1);
					\draw (w2)--(x2)--(y2);
					
					\node  at (1.5,1.5){\red{\xmark}}; 
					\node [rotate= 90]  at (2.2,0.75){\red{\xmark}}; 
					\end{tikzpicture}
					\subcaption{$(G_1, P_1) \oplus  (G_2, P_2)$}
				\end{minipage}
				\begin{minipage}[t]{0.33\textwidth}
					\centering
					\begin{tikzpicture}[>=latex,
					greensquare/.style={rectangle, draw=green, minimum size=1.5mm, inner sep=0pt},
					roundnode/.style={circle, draw=black,fill=black, very thick, minimum size=0.8mm, inner sep=0pt}]  
					\node [roundnode] (w1) at (0,1.5){};
					\node [roundnode] (x1) at (1,1.5){};
					\node [roundnode] (y1) at (2,1.5){};	
					\node [roundnode] (z1) at (2,0){};
					
					\node [roundnode] (w2) at (3.2,0){};
					\node [roundnode] (x2) at (2.2,0){};
					\node [roundnode] (y2) at (2.2,1.5){};	
					\node [roundnode] (z2) at (3.2,1.5){};
					\node [roundnode] (u) at (2,2.3){};
					
					\node [greensquare] (w1') at (0,1.5){};
					\node [greensquare] (x1') at (1,1.5){};
					\node [greensquare] (y1') at (2,1.5){};	
					\node [greensquare] (z1') at (2,0){};
					
					\node [greensquare] (y2') at (2.2,1.5){};	
					\node [greensquare] (z2') at (3.2,1.5){};
					
					\node at (0.2,1.25){$w_1$};
					\node at (1,1.25){$x_1$};
					\node at (1.8,1.25){$y_1$};
					\node at (1.9, -0.2){$z_1$};
					\node at (1, 0.75){$G_1$};
					\node at (2.2, 2.4){$u$};
					
					\node at (3.5,0){$w_2$};
					\node at (2.4,-0.25){$x_2$};
					\node at (2.3,1.75){$y_2$};
					\node at (3.4, 1.25){$z_2$};
					\node at (3, 0.75){$G_2$};
					
					\draw (x1)--(y1)--(z1);
					\draw (w2)--(x2)--(y2);
					\draw (x1)--(u);
					
					\draw [->] (w1)--(x1);
					\draw [->] (z2)--(y2);
					\draw [->] (2.5,0.5)--(x2);
					\draw [->] (2.7,0.25)--(x2);
					\draw [->] (2.5,0.9)--(y2);
					\draw [->] (2.7,1.2)--(y2);
					
					\draw [magenta,->] (u)--(z2);
					\draw [magenta, ->] (y1) to [bend right = 30] (x1);
					
					\node  at (1.5,1.45){\red{\xmark}}; 
					\node [rotate= 90]  at (2.2,0.75){\red{\xmark}}; 
					\node [rotate= 35] at (1.6,2){\textcolor{magenta}{\xmark}}; 
					\end{tikzpicture}
					\subcaption{$(G_1, P_1) \hat{\oplus} (G_2, P_2)$}
				\end{minipage}
				\begin{minipage}[t]{0.3\textwidth}
					\centering
					\begin{tikzpicture}[>=latex,	
					greensquare/.style={rectangle, draw=green, minimum size=1.5mm, inner sep=0pt},
					roundnode/.style={circle, draw=black,fill=black, very thick, minimum size=0.8mm, inner sep=0pt}]  
					\node [roundnode] (w1) at (0,1.5){};
					\node [roundnode] (x1) at (1,1.5){};
					\node [roundnode] (y1) at (2,1.5){};	
					\node [roundnode] (z1) at (2,0){};
					
					\node [roundnode] (u) at (1.5,2.3){};
					\node [roundnode] (v) at (2.5,2.3){};
					
					\node [roundnode] (w2) at (3.2,0){};
					\node [roundnode] (x2) at (2.2,0){};
					\node [roundnode] (y2) at (2.2,1.5){};	
					\node [roundnode] (z2) at (3.2,1.5){};
					
					\node [greensquare] (w1') at (0,1.5){};
					\node [greensquare] (x1') at (1,1.5){};
					\node [greensquare] (y1') at (2,1.5){};	
					\node [greensquare] (z1') at (2,0){};
					
					\node [greensquare] (y2') at (2.2,1.5){};	
					\node [greensquare] (z2') at (3.2,1.5){};
					
					\node at (0.2,1.25){$w_1$};
					\node at (1,1.25){$x_1$};
					\node at (1.8,1.25){$y_1$};
					\node at (1.9, -0.2){$z_1$};
					\node at (1, 0.75){$G_1$};
					
					\node at (1.3, 2.35){$u$};
					\node at (2.7, 2.35){$v$};
					
					\node at (3.5,0){$w_2$};
					\node at (2.4,-0.25){$x_2$};
					\node at (2.3,1.75){$y_2$};
					\node at (3.4, 1.25){$z_2$};
					\node at (3, 0.75){$G_2$};
					
					\draw (x1)--(y1)--(z1);
					\draw (w2)--(x2)--(y2);
					\draw [->] (z2)--(y2);
					\draw (u)--(v);
					
					\draw [->] (w1)--(x1);
					\draw [magenta, ->] (y1) to [bend right = 30] (x1);
					\draw [->] (2.5,0.9)--(y2);
					\draw [->] (2.7,1.2)--(y2);
					\draw [->] (2.5,0.5)--(x2);
					\draw [->] (2.7,0.25)--(x2);
					
					\node  at (2,2.3){\textcolor{magenta}{\xmark}};
					\draw [magenta, ->] (x1) to  (u);
					\draw [magenta, ->] (v) to  (z2);
					\node  at (1.5,1.5){\red{\xmark}}; 
					\node [rotate= 90]  at (2.2,0.75){\red{\xmark}}; 
					\end{tikzpicture}
					\subcaption{$(G_1, P_1) \tilde{\oplus} (G_2, P_2)$}
				\end{minipage}
				\caption{Illustration of the proof for \ref{formula-1}.}
				\label{fig-illustrations}
			\end{figure} 
			
			For \ref{formula-3}, let $M = M_1 \cup M_2$. If $\circ = \hat{\oplus}$, then let $D = D_1 \cup D_2 -(z_2,y_1) +(y_1,x_1)+ (y_1,z_2)+(z_2,u)$; If $\circ = \tilde{\oplus}$, then let $D = D_1 \cup D_2 - (z_2,y_1) +(y_1,x_1)+ (x_1,u) + (y_1,z_2)+(z_2,v)$.
			
			For \ref{formula-4}, let $M = M_1 \cup M_2$. If $\circ = \hat{\oplus}$, then let $D = D_1 \cup D_2 +(y_1,x_1)+ (z_2,u)$; If $\circ = \tilde{\oplus}$, then let $D = D_1 \cup D_2 + (y_1,x_1)+ (x_1,u) +(z_2,v)$.
			
			For \ref{formula-5}, let $D = D_1 \cup D_2 - (z_2,y_1)+ (y_1,x_1)+ (y_1,z_2) + (z_2,u)$ and $M = M_1 \cup M_2$.
			
			For \ref{formula-6}, let $M = M_1 \cup M_2$ and $D = D_1 \cup D_2 -(z_2,y_1) +(y_1,x_1)+ (y_1,z_2)+(z_2,u)$.
			
			In each of the cases, it is straightforward to verify that $(D,M)$ is an $(\vec{a}'', \vec{b}'')$-decomposition of $(G_1, w_1x_1y_1z_1) \circ (G_2, w_2x_2y_2z_2)$. 
		\end{proof} 
			
		  Instead of proving Lemma~\ref{main-lemmaA} directly, we  prove    Lemma \ref{lemma-property}, which is   more technnical, but facilitate  the usage of induction. Lemma \ref{main-lemmaA} will  follow from Lemma \ref{lemma-property},  except that the second half of Lemma~\ref{main-lemmaA}\ref{lemP} will be proven in the end of this section. Recall that 
			$$\tilde{\mathcal{G}} = \mathcal{G}(1101, 0000) \cap \mathcal{G}(1011, 0000) \cap \mathcal{G}(1002, 0000) \cap \mathcal{G}(2001, 0000).$$
			\begin{lemma}
				\label{lemma-property} 
				For the families $\fr$- and $\fp$-configurations, the following hold:
				\begin{enumerate}[label= {(\roman*)}] 
					\item \label{pro-R}
					$\fr\setminus \fr_2 \subseteq  \tilde{\mathcal{G}}   \cap \mathcal{G}(1001,1001)|_{zz^+\in M}  \cap  \mathcal{G}(1001,1100) \cap  \mathcal{G}(1001, 0011)|_{yz \in M}$. 
					\item \label{pro-R2} $\fr_3 \subseteq  \mathcal{G}(1001,0001)|_{zz^+\in M}$.  
					\item \label{pro-Q} $\fr_2 \subseteq   \mathcal{G}(1001,1001)|_{zz^+\in M}  \cap \mathcal{G}(1101,0000)  \cap \mathcal{G}(1011, 0000)|_{(z,y) \in D}  \cap \mathcal{G}(1001,1100)  \cap \mathcal{G}\\(1001, 0011)|_{(z,y) \in D}$.  
					\item \label{pro-P}	 $\fp \setminus \fp_2 \subseteq   \tilde{\mathcal{G}} \cap  \mathcal{G}(1001,0001)|_{zz^+ \in M} \cap \mathcal{G}(1001,1000)|_{w^-w \in M}$.
				\end{enumerate}
			\end{lemma}
			\begin{proof}
				The proof is by induction on $|V(G)|$. The case $(G,wxyz) \in \fr_1 \cup \fp_1$ is trivial. Assume  $(G,wxyz) = (G_1,w_1x_1y_1z_1) \circ (G_2, w_2x_2y_2z_2)$ for some $\circ \in \{\oplus, \hat{\oplus}, \tilde{\oplus}\}$. 
				
				Most of the statements directly follow from Lemma \ref{lem-operation} by induction on $(G_1,w_1x_1y_1z_1)$ and $(G_2, w_2x_2y_2z_2)$. More precisely,  
				\begin{itemize}
					\item  Lemma \ref{lem-operation} \ref{formula-0} implies that $\fr_2 \subseteq  \mathcal{G}(1001,1001)|_{zz^+\in M}$. 
					\item  Lemma \ref{lem-operation}\ref{formula-3} implies that $\fr_3, \fp_4,\fp_3 \subseteq \mathcal{G}(1001,0001)|_{zz^+ \in M}$. As a corollary, $\fp_4 \subseteq \mathcal{G}(1001,1000)|_{w'w \in M}$ (by the symmetry of $\fp_4$, see Observation \ref{obs-key}).
					\item  Lemma \ref{lem-operation}\ref{formula-1} implies that $\fr_3,\fp_4,\fp_3 \subseteq  \mathcal{G}(1011,0000)$ and $\fr_2 \subseteq   \mathcal{G}(1011,0000)$. Note that in the proof of Lemma \ref{lem-operation}\ref{formula-1}, we actually used the condition  that $(G_2,w_2x_2y_2z_2) \in \mathcal{G}(1001,1100)$, this implies that $(z,y) \in D$ in last formula. Therefore, $\fr_2 \subseteq   \mathcal{G}(1011,0000)|_{(z,y) \in D}$. We also observe that the last formula implies that $\fr_3  \subseteq  \mathcal{G}(1001,0011)|_{yz \in M}$ by Observation \ref{obs-RP and Q}. Also by  Observation \ref{obs-key}, $\fp_4 \subseteq  \mathcal{G}(1101,0000)$. 
					\item  Lemma \ref{lem-operation}\ref{formula-2} implies that $\fr_2 \subseteq   \mathcal{G}(1101,0000)$.
					\item Lemma \ref{lem-operation}\ref{formula-4} implies that $\fr_3, \fp_4,\fp_3 \subseteq  \mathcal{G}(1002,0000)$,  hence by Observation \ref{obs-key}, $\fp_4 \subseteq  \mathcal{G}(2001,0000)$.
					\item  Lemma \ref{lem-operation}\ref{formula-5} implies that $\fr_3, \fp_3 \subseteq  \mathcal{G}(2001,0000)$. 
					\item  Lemma \ref{lem-operation}\ref{formula-6} implies that $\fr_3, \fp_3 \subseteq  \mathcal{G}(1001,1000)|_{ww' \in M}$.  
				\end{itemize}
			
			 Now we show that $\fr_2 \subseteq  \mathcal{G}(1001,1100)$. By Observation \ref{obs-key}, $(G,z^+zyx) \in \fr_2$. As we have proved that $\fr_2 \subseteq  \mathcal{G}(1001,1001)$,  $(G,z^+zyx) $ has  a  $(1001,1001)$-decomposition $(D,M)$.    Then $(D-(x,y)+(z,y),M)$ is a  $(1001,1100)$-decomposition of $(G,wxyz)$.
				
			 Similarly, we can show that $\fr_3 \subseteq  \mathcal{G}(1001,1100)$ by using $\fr_3 \in \mathcal{G}(1001,0001)|_{zz^+ \in M}$ (Note that $(G,wxyz) \in \fr$ implies that  $(G,z^+zyx) \in \fr$ by Observation \ref{obs-key}).  
				
				Next we show that $\fr_2 \subseteq  \mathcal{G}(1001,0011)|_{(z,y) \in D}$. Assume $(G, wxyz) \in \fr_2$. Then  $(G, wxyz) = (G_1,w_1x_1y_1z_1)  \oplus (G_2, w_2x_2y_2z_2)$ for some   $(G_1,w_1x_1y_1z_1) \in \fr$ and $(G_2, w_2x_2y_2z_2) \in \fr$. By induction hypothesis, $(G_2, w_2x_2y_2z_2) \in \mathcal{G}(1001,1001)$.  Let $(D_2,M_2)$ be a $(1001,1001)$-decomposition of $(G_2, w_2x_2y_2z_2)$. 
				Depending on $(G_1,w_1x_1y_1z_1) \in \fr_1\cup \fr_3$ or $(G_1,w_1x_1y_1z_1) \in   \fr_2$, we consider two cases.

			\begin{itemize}
				\item [-]	If $(G_1,w_1x_1y_1z_1) \in \fr_1 \cup \fr_3$, then  $(G_1,w_1x_1y_1z_1) \in \mathcal{G}(1001,0011)|_{y_1z_1 \in M}$   by induction hypothesis. Assume $(D_1,M_1)$  is a $(1001,0011)$-decomposition of $(G_1, w_1x_1y_1z_1)$ with  $y_1z_1\in M_1$. Then $(D_1 \cup D_2,M_1 \cup M_2)$ is a $(1001,0011)$-decomposition of $(G,wxyz)$ with $(z,y) \in D$. 
				
				\item [-]	If  $(G_1,w_1x_1y_1z_1) \in \fr_2$, then $(G_1,w_1x_1y_1z_1) \in \mathcal{G}(1001,0011)|_{(z_1,y_1) \in D}$ by induction hypothesis. Assume $(D_1,M_1)$  is a $(1001,0011)$-decomposition of $(G_1, w_1x_1y_1z_1)$ with $(z_1,y_1) \in D_1$. Then $(D_1 \cup D_2,M_1 \cup M_2)$ is also a  $(1001,0011)$-decomposition of $(G,wxyz)$ with $(z,y) \in D$.
			\end{itemize}

				Next we show that $\fr_3 \subseteq  \mathcal{G}(1101,0000)$.
				Assume $(G, wxyz) \in \fr_3$. By Observation \ref{obs-associative},   $(G, wxyz) = (G_1,w_1x_1y_1z_1)  \hat{\oplus} (G_2, w_2x_2y_2z_2)$ for some   $(G_1,w_1x_1y_1z_1) \in \fr$ and $(G_2, w_2x_2y_2z_2) \in \fr$ (note that here $w=w_1$, $x=x_1$, $y=u$ and $z=z_2$).
				Depending on $(G_1,w_1x_1y_1z_1) \in \fr_1 \cup \fr_3$ or $(G_1,w_1x_1y_1z_1) \in  \fr_2$, we consider two cases.
				\begin{itemize}
					\item[-] If $(G_1, w_1x_1y_1z_1) \in \fr_1 \cup \fr_3$, then by induction hypothesis,    both $(G_1, w_1x_1y_1z_1)$ and $(G_2, w_2x_2y_2z_2)$ have $(1101,0000)$-decompositions, say $(D_1,M_1)$ and $(D_2,M_2)$, respectively.  Let $D= D_1\cup D_2-(z_1,y_1)-(z_2,y_1)+(y_1,x_1)+(y_1,z_2)+(z_2,u)$, $M= M_1\cup M_2 + z_1y_1$. Note that $y_2$ ($= y_1$) is a source in $D$ (as $d_G(y_1) = d_{G_1}(y_1)+d_{G_2}(y_2)-1=3$ and $(y_1,z_2), (y_1,x_1) \in D$ and $y_1z_1 \in M$). So $D$ is acyclic, and hence $(D,M)$ is a $(1101,0000)$-decomposition of $(G,wxyz)$. 
					
					\item[-] If $(G_1, w_1x_1y_1z_1) \in \fr_2$, then by induction hypothesis,   $(G_1, w_1x_1y_1z_1) \in\mathcal{G}(1001,0011)|_{{(z_1,y_1) \in D}}$  and $(G_2, w_2x_2y_2z_2) \in \mathcal{G}(1101,0000)$. Let $(D_1,M_1)$ be a $(1001,0011)$-decomposition of $(G_1, w_1x_1y_1z_1)$ with  {$(z_1,y_1) \in D_1$,} and $(D_2,M_2)$ be a  $(1101,0000)$-decomposition of $(G_2, w_2x_2y_2z_2)$. Let $D= D_1\cup D_2-(z_1,y_1)-(z_2,y_2)+(x_1,y_1)+(y_1,z_1)+(y_1,z_2)+(z_2,u)$, $M= M_1\cup M_2$. Note that in $D$ restricted to $G_1$, $y_1$ has out-degree $1$, and $z_1$ is the head of $(y_1,z_1)$, but $z_1$ has out-degree $0$ in $G_1$, so there is no directed cycle in $D$ restricted to $G_1$. On the other hand, $y_1$ is a source of $D$ restricted to $G_2$, and $z_2$ has out-degree $0$ in $D$ restricted to $G_2$. So $D$ is acyclic, hence $(D,M)$ is a $(1101,0000)$-decomposition of $(G,wxyz)$.
				\end{itemize}
		 
				Last we show that $\fp_3 \subseteq  \mathcal{G}(1101,0000)$. Assume  $(G, wxyz) \in \fp_3$. Then $(G, wxyz) = (G_1,w_1x_1y_1z_1)  \hat{\oplus} (G_2, w_2x_2y_2z_2)$ for some   $(G_1,w_1x_1y_1z_1) \in \fp$ and $(G_2, w_2x_2y_2z_2) \in \fr$. By induction hypothesis, $(G_1,w_1x_1y_1z_1) \in  \mathcal{G}(1001,0001)$. By Observation \ref{obs-key} and induction hypothesis, $(G_2, z^+_2z_2y_2x_2) \in \fr \subseteq   \mathcal{G}(1011,0000)$. Let $(D_1,M_1)$ be a $(1001,0001)$-decomposition of $(G_1, w_1x_1y_1z_1)$, and let $(D_2,M_2)$ be a $(1011,0000)$-decomposition of $(G_2, z^+_2z_2y_2x_2)$. Let $D= (D_1-(z_1,y_1))\cup D_2 + (x_1,y_1)+(y_1,z_2)+(z_2,u)$,   and {$M= M_1\cup M_2$.} As  $d_{D_2}(z_2)=0$,	$d_{D_1}(z_1)=1$ and $(z_1,y_1) \in D_1$, we know that $D$ is acyclic. Therefore, $(D,M)$ is a $(1101,0000)$-decomposition of $(G,wxyz)$. 
			\end{proof}
			
	 Now we prove the second half Lemma~\ref{main-lemmaA}\ref{lemP}.     Assume   $(G,wxyz) \in \fp_4$. we need to show that $(G,w^-wxy) \in   \mathcal{G}(1001,0000)$. 
	 By Observation~\ref{obs-RP and Q}, $(G-\{x,y\},w^-wuz) \in \fr_2$, where $u$ is the vertex in $V(G)\setminus V(B_G)$, and $uwxyz$ is a $5$-face of $G$. By Lemma \ref{lemma-property}\ref{pro-Q}, $(G-\{x,y\},w^-wuz)$ has a 
	 $(1001,0011)$-decomposition $(D,M)$ with $(z,u) \in D$. Let $D'=D-(z,u)+(u,z)+(z,y)+(y,x)+(u,w)$ and $M'=M$. Then $(D',M')$ is a $(1001,0000)$-decomposition of $(G,w^-wuz)$.
	 
	 This completes the proof of Lemma~\ref{main-lemmaA}.

\end{document}